
\input amstex
\documentstyle{amsppt}
\magnification=1200
\NoBlackBoxes
\TagsOnRight
\refstyle {A}
\widestnumber\key{BDGW}
\def \-->{\longrightarrow}
\def \kler{K\"ahler}

\def \*bar{\overline{*}}

\def \curvdh{F_{\dbare,H}}

\def \dbar{\overline{\partial}}
\def \dbare{\overline{\partial}_E}

\def \G{\frak G}
\def \Gc{\frak G_\Bbb C}

\def \lieG{\frak g}

\def\astable{$\alpha$-stable}
\def\astability{$\alpha$-stability}

\def \E{\Cal E}
\def \C{\Cal C}

\def \cE{\Cal E}
\def \cExt{\Cal Ext}
\def \cF{\Cal F}

\def \bO{\bold O}
\def \extn{0\lra \cE_1\lra \cE\lra \cE_2\lra 0}
\def \sextn{0\lra \cE'_1\lra \cE'_1\lra \cE'_2\lra 0}
\def \fextn{0\lra \cF_1\lra \cF\lra \cF_2\lra 0}
\def \he{Hermitian-Einstein\ }
\def \lra{\longrightarrow}
\def \hra{\hookrightarrow}

\def \Met{\text{Met}}
\def \ct{cohomology triple}
\def \tr{\text{Tr}}
\def \log{\text{log}}
\def \vol{\text{vol}}
\def \tri{(\cE_1,\cE_2,\Phi)}
\def \stri{(\cE'_1,\cE'_2,\Phi')}
\def \ftri{(\cF_1,\cF_2,\Psi)}
\def \fstri{(\cF'_1,\cF'_2,\Psi')}
\def \hpt{H^p(Hom (\cE_2,\cE_1))}
\def \Hom{\text{Hom}}
\def \cHom{Hom}
\def \ext{\text{Ext}}
\def \ast{*}
\def \as{\astable}
\def \asty{\astability}
\def \*E{\overline{*}_E}
\def \im{\text{Im}}
\def \Par{\bold{Par}}
\def\today{\ifcase\month\or
  January\or February\or March\or April\or May\or June\or
  July\or August\or September\or October\or November\or
  December\fi
  \space\number\day, \number\year}

\topmatter
\title
Higher cohomology triples and holomorphic extensions
\endtitle
\author
 Steven B. Bradlow and Oscar Garcia-Prada
\endauthor
\address
Department of Mathematics, University of Illinois,
Urbana, IL 61801
\endaddress
\email  bradlow@uiuc.edu \endemail
\address Departamento de Matematicas,Universidad Autonoma de Madrid,\\ Ciudad
Universitaria de Cantoblanco, 28049 Madrid,SPAIN \endaddress
\email OGPRADA@ccuam3.sdi.uam.es\endemail
\abstract We introduce equations for special metrics, and notions of stability
for some new types of augmented holomorphic bundles. These new examples
include
holomorphic extensions, and in this case we prove a Hitchin-Kobayashi
correspondence between a certain deformation of the Hermitian-Einstein
equations and our definition of stability for an extension.
\endabstract
\endtopmatter
\document

\heading\S 1. Introduction\endheading

Let $E\-->X$\ be a fixed smooth bundle over a \kler\ manifold. There are three
natural moduli spaces associated to $E$; one algebraic, one complex analytic,
and one symplectic. The first is the moduli space of slope stable holomorphic
structures on $E$, and is constructed by Geometric Invariant Theory. The
second, the moduli space of Hermitian-Einstein connections, is constructed by
gauge theory and deformation theory. For the third, one uses the symplectic
structure induced on the space of unitary connections and considers the moment
map for the action of the unitary gauge group. The symplectic moduli space is
then the Marsden-Weinstein quotient of the zero level of the moment map by the
action of the unitary gauge group.

In fact, it is well known that these three quotients can all be identified,
and
this is referred to as the Hitchin-Kobayashi correspondence.  Furthermore,
this
triad of descriptions has, in  recent years, been found to be a common feature
in an ever expanding range of situations.
In most of these, the moduli spaces are for augmented
bundles of one kind or another, i.e. for objects consisting
of one or more holomorphic bundle together with
prescribed holomorphic sections. A summary of such results can be found in
[BDGW].

In this paper we discuss some extensions of these ideas
in
two directions that have not hitherto been pursued.  This
involves consideration of an interesting class of
equations
which includes deformations of the Hermitian-Einstein
equations as well as certain generalizations of the equations known as the
vortex equations. It also requires the introduction of a new  concepts of
stability for various augmented bundles.

In the one class of examples that we discuss, the starting
point is the observation that symplectic reduction can be
carried out more generally than simply at the 0-level
set. In particular, symplectic quotients can be constructed from the inverse
images of coadjoint orbits in the dual of
the Lie algebra of the unitary gauge group. It is natural to look for a
description of such reduced spaces as complex quotients and to try to find an
algebraic characterization of this quotient as a moduli space.

The simplest example of such a generalization can be
described as follows.  Suppose that $E=E_1\oplus E_2$.
Fix a smooth metric $K$\ on $E$\ such that the above splitting
of
$E$\ is an orthogonal decomposition. Let
$T_{\tau_1,\tau_2}\in\G$\ be the global gauge
transformation
given by
$$T_{\tau_1,\tau_2}= \pmatrix
i\tau_1 \bold I_1 & 0 \\
0 & i\tau_2 \bold I_2
\endpmatrix\ $$
with respect to the given splitting of $E$, and let $\Cal
O(\tau_1,\tau_2)$\ be the coadjoint orbit of
$T_{\tau_1,\tau_2}$. The points in the
inverse image $\Psi^{-1}(\Cal O(\tau_1,\tau_2))$\ can be
described as \it holomorphic extensions\rm, either of
the
type
$$0\-->\Cal E_1\-->\Cal E\-->\Cal E_2\-->0\ ,$$
or of the type
$$0\-->\Cal E_2\-->\Cal E\-->\Cal E_1\-->0\ .$$
Pursuing this example, we find interesting
``deformations"
of both the Hermitian-Einstein equations and the notion
of
bundles stability.  Furthermore, these are naturally
interpreted in terms of holomorphic extensions.

The second type of structure we consider is a natural generalization of the
triples described in [BGP]. In [BGP]
we described objects consisting of two holomorphic bundles,
$\Cal E_1$\ and $\Cal E_2$\ plus a map between them, i.e. a section $\Phi\in
H^0(X,Hom(\Cal E_2,\Cal E_1))$.  In the generalization we have in mind, we
take
$\Phi$\  in $H^p(X,Hom(\Cal E_2,\Cal E_1))$, for any $p$.
We call such objects \it p-cohomology
triples\rm.  Apart from their interest as natural generalizations of the
original triples, such objects (with $p=2$) have been encountered in the work
of Pidstrigach and Tyurin ([PT]), and more recently in connection with the
Seiberg-Witten invariants for algebraic surfaces (cf.[W]).  For the case
$p=0$,
we described in [BGP] what the natural notion of stability is, and what the
corresponding equations for special metrics look like.  In this paper we
discuss  how these can be modified to describe the more general situation.

The case of $p=1$\ is of particular interest, since elements in
$H^1(X,Hom(\Cal
E_2,\Cal E_1))$\ can be interpreted as extension classes.  This leads to
interesting relations between the two kinds of situations described above. We
describe in some detail how these points of view compare.  We also relate
these
to yet another description of holomorphic extensions, namely one in terms of
the bundles $\cE$, $\cE_2$\ plus surjective maps $\pi:\cE\-->\cE_2$.  Such
objects, which describe extensions of $\cE_2$\ by the kernel of the map, can
be
thought of as a special type of $p=0$\ triples. More specifically, they
correspond to such triples in which the map between the bundles is surjective.
We thus discuss the relation between such surjective ($p=0$) triples,
1-cohomology triples and extensions.

\it Remark.\rm\   The result given in Theorem 3.9 has been proved
independently
by Daskalopoulos, Uhlenbeck and Wentworth [DUW]. With stability defined as in
in Definition 3.4, they have gone on to give analytic as well as invariant
theory constructions of the moduli spaces of stable extensions.

\subheading{Acknowledgements} The authors would like to
thank Alastair King for many helpful conversations and ideas, especially with
regard to the formulation of the definitions of stability in \S 2.2.  Both
authors are members of the VBAC group of Europroj.

\heading \S 2 Cohomology Triples\endheading

Let $(X,\omega)$ be a compact \kler\ manifold of
dimension $n$, and fix two smooth complex bundles $E_i\-->X,\ i=1,2$.  Denote
their ranks and degrees by $d_i$\ and $r_i$, where by the degree
we mean, in general, $\int_X{c_1(E_1)\wedge\omega^{n-1}}$.
In order to simplify certain formulae, we assume that the volume of $X$\ is
normalized to $2\pi$.
A holomorphic triple based on $E_1$\ and $E_2$\ consists of holomorphic
structures (given by $\dbar$-operators $\dbar_1$\ and $\dbar_2$) on these
bundles plus
a holomorphic section of $Hom(E_2,E_1)$. Thus the
augmentation is represented by a smooth section
$\Phi\in\Omega^0(Hom(E_2,E_1))$\ satisfying the constraint
$\dbar_{1,2}(\Phi):=\dbar_1\circ\Phi-\Phi\circ\dbar_2=0$.

There are two distinct ways in which one might want to
generalize this to allow form-valued augmentations.
\roster
\item In the first, which we call
\it p-cocycle triples\rm, one replaces holomorphic
sections
of $\Omega^0(Hom(E_2,E_1))$\ by \it holomorphic
sections of
$\Omega^{0,p}(Hom(E_2,E_1)$\rm.
\item In the second, which we call
\it p-cohomology triples\rm, the augmentation is
considered
to be the class in $H^{0,p}(Hom(E_2,E_1))$\ represented
by a
holomorphic section in $\Omega^{0,p}(Hom(E_2,E_1))$.
\endroster

As will be seen (cf. Section 2.2 ), there are compelling
reasons for regarding the second approach as the
``correct" one. Nevertheless, at least in the case where $X$\ is a  Riemann
surface, there are interesting features of both types of augmentation. In the
case that $\E_2$\ is  fixed to be the structure  sheaf, the resulting objects
may be considered as p-cocycle- and p-cohomology pairs.

\subheading{\hfil\S 2.1 The basics}

Set
$$\chi^{(p)}
=\C_1\times\C_2\times\Omega^{0,p}(Hom(E_2,E_1)\
,\tag 2.1.1$$
where $\C_i$\ denotes the space of holomorphic
structures
(or equivalently, the space
of $\dbar$-operators) on $E_i$.

\proclaim{Definition 2.1} We can define the
the space of all \bf p-cocycle triples\it\ on
$(E_1,E_2)$\
by the holomorphic subspace
$$\Cal Z^{(p)}=\{(\dbar_1,\dbar_2,\phi)\in\chi^{(p)}\ :\
\dbar_{1,2}(\phi)=0\}\ ,\tag 2.1.2$$
where $\dbar_{1,2}(\phi)=\dbar_1\circ\phi-
\phi\circ\dbar_2$.
\endproclaim

We can define an equivalence relation on $\Cal Z^{(p)}$\
by
$$(\dbar_1,\dbar_2,\phi)\sim(\dbar_1,\dbar_2,\phi+
\dbar_{1,2}(\alpha))\ ,$$
for any $\alpha\in \Omega^{0,p-1}(Hom(E_2,E_1)$.  The
p-cohomology triples
are described by the equivalence classes in  $\Cal
Z^{(p)}/\sim$.  Notice that these equivalence
classes correspond to orbits
of the additive group $\Omega^{p-1}_{1,2}:=
\Omega^{0,p-1}(Hom(E_2,E_1)$ under the action
$$\alpha\circ(\dbar_1,\dbar_2,\phi)=
(\dbar_1,\dbar_2,\phi+\dbar_{1,2}(\alpha))\ .$$

\proclaim{Definition 2.2}  The space of all \bf p-
cohomology
triples\it\ on $(E_1,E_2)$\ is defined by

$$\Cal H^{(p)}=\Cal Z^{(p)}/\Omega^{p-1}_{1,2}\ .\tag 2.1.3$$
\endproclaim

\proclaim{Definition/Lemma 2.3}The complex gauge
group
$\Gc=\Gc^{(1)}\times\Gc^{(2)}$\
acts on both $\Cal Z^{(p)}$\ and $\Cal H^{(p)}$. In
both cases, the $\Gc$-orbits
correspond to isomorphism classes, with the notion of
isomorphism defined in the obvious way. Thus the
``moduli
spaces" of isomorphism classes of p-cocycle( resp.
p-cohomology) triples corresponds to the orbit
space $\Cal Z^{(p)}/\Gc$( resp.
$\Cal H^{(p)}/\Gc$).
\endproclaim

It is important to observe that in the double quotient
$$\Cal H^{(p)}/\Gc=(\Cal Z^{(p)}/
\Omega^{p-1}_{1,2})/\Gc\ ,$$
the order of the quotient operations
cannot be reversed.  Not only do the actions of
$\Omega^{p-1}_{1,2}$\ and
$\Gc$\ fail to commute, but
$\Omega^{p-1}_{1,2}$\ does not act on
$\Cal Z^{(p)}/\Gc$\ in any obvious way. Nevertheless,
the quotient $\Cal H^{(p)}/\Gc$\ can be described
as a quotient of $\Cal Z^{(p)}$, namely as the quotient
by the group action of the semidirect product
$\Omega^{p-1}_{1,2}\ltimes\Gc$.

\proclaim{Definition/Lemma 2.4} We can identify
$$(\Cal Z^{(p)}/\Omega^{p-1}_{1,2})/\Gc= \Cal
H^{(p)}/\Omega^{p-1}_{1,2}\ltimes\Gc\ $$
where the group structure on the semidirect product is
defined by
$$(\alpha,g_1,g_2)(\alpha',g'_1,g'_2)=
({g'_1}^{-1}\alpha g'_2+\alpha',g_1g'_1,g_2g'_2)\ ,\tag 2.1.4$$
and the action on $\Cal Z^{(p)}$\ is
$$ (\alpha,g_1,g_2)(\dbar_1,\dbar_2,\phi)=
(g_1(\dbar_1),g_2(\dbar_2),
g_1(\phi+\dbar_{1,2}(\alpha))g_2^{-1})\ .\tag 2.1.5$$
\endproclaim

\subheading{\hfil\S 2.2 Stability with parameters}

As usual, one cannot expect the orbit spaces
$\Cal Z^{(p)}/\Gc$\ or
$\Cal H^{(p)}/\Gc$\ to yield well behaved moduli
spaces without
restricting to suitably defined spaces of ``stable" orbits.
The definition of stability that we propose for
cohomology-triples is a reasonably
straightforward extensions of the stability defined for
triples in [BGP].
Since this definition is in terms of a condition on
subtriples, we need to
specify precisely what we mean by the subobjects of
cohomology  triples.

Let $\cE_1=(E_1,\dbar_1)$ and $\cE_2=(E_2,\dbar_2)$
holomorphic vector bundles
on $X$ and $\Phi\in \hpt$. To define the subobjects of
the
\ct\ $T=\tri$ we need
to determine the category to which $T$ belongs. The
subobjects of $T$ will be
then certain objects in this category---the ones for
which
there is an injective morphism
to $T$.

The category we need to consider is the category of
``$\ext^p$'' triples. Its elements
consist of triples $(\cF_1,\cF_2,\Psi)$, where $\cF_1$
and
$\cF_2$ are coherent sheaves
on $X$ and  $\Psi$ is an element of
$\ext^p(\cF_2,\cF_1)$.

Recall that
$\ext^0(\cF_2,\cF_1)=\Hom(\cF_2,\cF_1)\cong
H^0(\cHom(\cF_2,\cF_1))$, and
if $\cF_2$ is locally free:

1) $\cHom(\cF_2,\cF_1)\cong \cF_1\otimes\cF_2^\ast$

2) $\ext^p(\cF_2,\cF_1)\cong
H^p(\cF_1\otimes\cF_2^\ast)$.

Let $T=(\cF_1,\cF_2,\Phi)$ and
$T'=(\cF_1',\cF_2',\Psi')$ be
two $\ext^p$ triples.
A morphism $T'\lra T$ consists of morphisms
$f_1:\cF_1'\lra\cF_1$ and
$f_2:\cF_2'\lra\cF_2$ such that under the induced maps

$$\CD
\ext^p(\cF_2', \cF_1') @>{f_{1\ast}}>>
\ext^p( \cF_2', \cF_1)
@<{f_2^\ast}<< \ext^p( \cF_2, \cF_1)
\endCD$$

One has that
$$
f_{1\ast}(\Phi')=f_2^\ast (\Phi).\tag 2.2.1
$$

Note that when $p=0$,  this is equivalent to having the
following commutative
diagram

 $$\CD
 \cF_2@>{\Phi}>> \cF_1 \\
 @AA{f_2}A  @AA{f_1} A\\
 \cF_2'@>{\Phi'}>>\cF_1'.
 \endCD
 $$

\proclaim{Definition 2.5}Let $T=\ftri$ be an $\ext^p$ triple. A subobject $T'$
of
$T$
consists of
an $\ext^p$ triple $\fstri$ such that one has injections
$i_1:\cF_1'\hra\cF_1$ and
$i_2:\cF_2'\hra\cF_2$, which induce a morphism from
$T$ to
$T'$, i.e.
$i_{1\ast}(\Psi')=i_2^\ast(\Psi)$.
\endproclaim

Notice that if $\tri$ is a \ct, (where
$\cE_1$ and $\cE_2$ are locally
free and $\Phi\in\hpt$), a subobject NEED NOT BE  a
cohomology triple. In general it will be only an $\ext^p$
triple, i.e. $\cE_1'\hra \cE_1$ and $\cE_2'\hra\cE_2$ are
not necessarily locally free
and $\Phi'\in \ext^p(\cE_2',\cE_1')$ (If $\cE_2'$ is not
locally free this cannot be
identified with $H^p(\cHom(\cE_2',\cE_1'))$.

The definition of stability is given in terms of \it defect functions\rm\
(Alastair King's terminology)
for pairs of bundles:

\proclaim{Definition 2.6} Let $(E_1,E_2)$\ be a pair
of bundles of degree $d_1$ and $d_2$, and rank $r_1$\
and
$r_2$.
Fix real numbers $\{a_1,a_2,\tau_1,\tau_2\}$, and
define
$\theta_{a_1,a_2,\tau_1,\tau_2}(E_1,E_2)$\ by
$$\theta_{a_1,a_2,\tau_1,\tau_2}(\E_1,\E_2)=
a_1d_1+a_2d_2-\tau_1r_1-\tau_2r_2\ .\tag 2.2.2$$
\endproclaim

\proclaim{Definition 2.7} Let
$(\E_1,\E_2,\Phi)$
be a  p-cohomology triple based on the
smooth bundles
$(E_1,E_2)$. Fix real numbers
$\{a_1,a_2,\tau_1,\tau_2\}$\ with $a_1$\ and $a_2$\ non-negative, and such
that
$$a_1d_1+a_2d_2-\tau_1r_1-\tau_2r_2=0\ ,$$
i.e. such that
$\theta_{a_1,a_2,\tau_1,\tau_2}(\E_1,\E_2)=
0$. We
say that the triple
$(\E_1,\E_2,\Phi)$\ is \bf $\{a_1,a_2,\tau_1,\tau_2\}$-
stable\it\ if
$$\theta_{a_1,a_2,\tau_1,\tau_2}(\E'_1,\E'_2)< 0\ $$
for all p-cohomology subtriples.
\endproclaim

\noindent\bf Remarks:\rm\ \roster\item As usual, to study stability questions
it suffices to consider \it
saturated\rm\ subobjects of
$\tri$ (we are assuming that $\cE_1$ and $\cE_2$ are
torsion
free). These are subobjects
$\stri$ for which the inclusions $\cE_1'\hra\cE_1$ and
$\cE_2'\hra\cE_2$ are saturated,
i.e. $\cE_1/\cE_1'$ and $\cE_2/\cE_2'$ are torsion free.

\item As in the case of \it cohomology\rm\ triples, before we
can define stability for a \it cocycle\rm\ triple, we must
first establish what the legitimate subobjects
are . We immediately run into
difficulty when we consider what the appropriate
category for such objects should be. Denote the
objects in the category by  $(\Cal F_1,\Cal F_2,\phi)$.
The problem is that we do not want to require that
$\Cal F_1$\ and $\Cal F_2$\ be locally free sheaves.
This means we need to have a replacement for
$\Omega^p(Hom(\Cal F_2,\Cal F_1))$\ in the case that
$\Cal F_1$\ and $\Cal F_2$\ are NOT locally free.  This
is
one of the reasons that cohomology triples are to be
preferred to their cocycle cousins .

Notice that when $X$\ is a Riemann surface, these
difficulties do not
arise, and sensible definitions can be given.  Subobjects
are defined to be cocycle triples $(\Cal F'_1,\Cal
F'_2,\phi')$
with injections $i:\Cal F'_1\-->\Cal F_1$\ and
$j:\Cal F'_2\-->\Cal F_2$\ such that
$$i\circ\phi'=\phi'\circ j\ .$$
Stability with respect to parameters
$\{a_1,a_2,\tau_1,\tau_2\}$\ is then defined exactly as
for
cohomology triples.

\item Finally, suppose that $(\E_1,\E_2,\phi)$\ is a $0$-cohomology
triple, i.e. a triple
 with $\phi\in H^0(Hom(E_2,E_1))$. We recover the old
definition of $\tau$-stability given in [BGP] by taking
$\{a_1,a_2,\tau_1,\tau_2\}=\{1,1,\tau,\tau'\}$.  The
definition
above is thus a generalization of $\tau$-stability.
\endroster

The  parameter space for the parameters in the definition of stability can be
described as follows.
Let $\Par\subset \Bbb R^4$\ be the subspace
$$\Par=\{(a_1,a_2,\tau_1,\tau_2)\ |\
a_1\ge 0\ ,\ a_2\ge 0\ ,\ a_1d_1+a_2d_2-\tau_1r_1-\tau_2r_2=0\ \}\ .\tag
2.2.3$$
Notice that the definition of $\{a_1,a_2,\tau_1,\tau_2\}$-
stability is insensitive to an overall scaling of $(a_1,a_2,\tau_1,\tau_2)$\
by
a positive scale factor.
The effective parameter space is thus $\Par/\Bbb R^+$.
The ``geography" of this parameter space is an interesting issue, which we
will
return to in a later paper. There are however, a few features which are
immediately apparent.

The first feature comes from the fact that (at least for the case when $X$\ is
algebraic), the degrees and ranks of subobjects may be assumed to be integers,
i.e. to lie in a discrete subset of $\Bbb R$. It follows immediately that

\proclaim {Lemma 2.8}
The parameter space $\Par/\Bbb R^+$\ is partitioned into chambers.  The walls
are determined by the choices of $(a_1,a_2,\tau_1,\tau_2)$\ at which the
relation $\theta_{a_1,a_2,\tau_1,\tau_2}(\E'_1,\E'_2)= 0$\ is numerically
possible.  Within a fixed chamber the definition of
$(a_1,a_2,\tau_1,\tau_2)$-stability is independent of the values of
$(a_1,a_2,\tau_1,\tau_2)$
\endproclaim

The next result identifies a special region within $\Par/\Bbb R^+$.

\proclaim{Proposition 2.9}
\roster
\item Suppose that $a_1>0$. Then the space of $(a_1,a_2,\tau_1,\tau_2)$-stable
objects is empty unless $\tau_1/a_1>\mu(\cE_1)$.
\item There are positive numbers $\epsilon_1,\epsilon_2$\ such that the
following is true:

Let $(a_1,a_2,\tau_1,\tau_2)$\ be any point in $\Par$\ such that
$$\align
&a_1a_2\ne  0\ ,\\
&0<\frac{\tau_1}{a_1}-\mu(E_1)<\epsilon_1\ ,\\
&\frac{\tau_1}{a_1}-\mu(E_1)<
 (\frac{a_2}{a_1})\epsilon_2\ .
\endalign$$
Then in any $(a_1,a_2,\tau_1,\tau_2)$-stable object, say
$(\cE_1,\cE_2,\Phi)$, the bundles $\cE_1$\ and $\cE_2$\ are semistable.
Conversely, if $\cE_1$\ and $\cE_2$\ are  stable bundles, then all cohomology
triples $(\cE_1,\cE_2,\Phi)$\ are $(a_1,a_2,\tau_1,\tau_2)$-stable.
\endroster
\endproclaim

\demo{Proof}  Both parts of the proposition use the following observations.
Let
$(\cE_1,\cE_2,\Phi)$ be any p-cohomology triple. For any subsheaf $\cE'_1$\ of
$\cE_1$\ we can construct the subtriple $(\cE'_1,0,0)$.  Furthermore, for any
subsheaf $\cE'_2$\ of $\cE_2$\ we can construct the subtriple
$(\cE_1,\cE'_2,\Phi')$\ by taking $\Phi'=i_2^*(\Phi)$.  Notice that
$$\theta_{a_1,a_2,\tau_1,\tau_2}(\E'_1,0)
=a_1r'_1(\mu(\cE'_1)-\tau_1/a_1)\ .\tag 2.2.4$$
Part (1) follows immediately from this. For part (2), we observe that if
$a_1\ne 0$, then we can write (2.2.4) as
$$\frac{\theta_{a_1,a_2,\tau_1,\tau_2}(\E'_1,0)}{a_1}
=r'_1(\mu(\cE'_1)-\mu(\cE_1))+
r'_1(\mu(\cE_1)-\tau_1/a_1)\ .\tag 2.2.5$$
It follows that if $|\mu(\cE_1)-\tau_1/a_1|$\ is sufficiently small, and if
$\theta_{a_1,a_2,\tau_1,\tau_2}(\E'_1,0)<0$, then $\mu(\cE'_1)-\mu(\cE_1)\le
0$.
We now consider the subobjects coming from subsheaves of $\cE_2$.  For these,
we get
$$\theta_{a_1,a_2,\tau_1,\tau_2}(\E_1,\cE'_2)
=a_1d_1+a_2d'_2-r_1\tau_1-r'_2\tau_2\ .\tag 2.2.6$$
Using the constraint equation $a_1d_1+a_2d_2-\tau_1r_1-\tau_2r_2=0$, this can
be written as
$$\frac{\theta_{a_1,a_2,\tau_1,\tau_2}(\E_1,\cE'_2)}{a_2}=
r'_2(\mu(\cE'_2)-\mu(\cE_2))-
r_1(\frac{r_2-r'_2}{r_2})
(\frac{\tau_1/a_1-\mu(\cE_1)}{a_2/a_1})\ .\tag 2.2.7$$
Thus if
$\theta_{a_1,a_2,\tau_1,\tau_2}(\E_1,\cE'_2)<0$, then
$$\mu(\cE'_2)-\mu(\cE_2)<(\frac{r_1}{r'_2})
(\frac{r_2-r'_2}{r_2})
(\frac{\tau_1/a_1-\mu(\cE_1)}{a_2/a_1})\ .$$
Since $(\frac{r_1}{r'_2})(\frac{r_2-r'_2}{r_2})$\ is bounded above, it follows
that $\mu(\cE'_2)-\mu(\cE_2)\le 0$\ if
$\frac{\tau_1/a_1-\mu(\cE_1)}{a_2/a_1}$\
is sufficiently small. This completes the proof of the first claim in (2).
The second claim also follows from the identities
(2.2.5) and (2.2.7), which show that for any subtriple $(\cE_1,\cE'_2,\Phi')$\
we have
$$\align
\theta_{a_1,a_2,\tau_1,\tau_2}(\E'_1,\cE'_2)&=
\theta_{a_1,a_2,\tau_1,\tau_2}(\E_1,\cE'_2)+
\theta_{a_1,a_2,\tau_1,\tau_2}(\E'_1,0)-a_1r_1(\mu(\cE_1)-\tau_1/a_1)\\
&=a_2r'_2(\mu(\cE'_2)-\mu(\cE_2))+
a_1r'_1(\mu(\cE'_1)-\mu(\cE_1))+\\
&+ a_1(r_1-r'_1)(\tau_1/a_1-\mu(\cE_1))
-a_2r_1(\frac{r_2-r'_2}{r_2})
(\frac{\tau_1/a_1-\mu(\cE_1)}{a_2/a_1})\ .
\endalign$$

\enddemo

\subheading {\hfil\S 2.3 Comparison of cocycle and
cohomology}

In the case that $X$\ is a Riemann surface, the following
comparison
between cohomology and cocycle triples makes sense.
Let $\pi:\Cal Z^{(p)}\-->\Cal H^{(p)}=\Omega^{p-
1}_{1,2}$\ denote the projection map.
Let $(a_1,a_2,\tau_1,\tau_2)$\ be any set of real
numbers
satisfying the constraint
$$a_1d_1+a_2d_2-\tau_1r_1-\tau_2r_2=0\ .$$

\proclaim{Proposition 2.10} Suppose that $X$\ is a
Riemann
surface.
For any cohomology triple
$(\dbar_1,\dbar_2,\Phi)\in\Cal H^{(p)}$, the following
are
equivalent:
\roster
\item  $(\dbar_1,\dbar_2,\Phi)\in\Cal H^{(p)}$\ is
$(a_1,a_2,\tau_1,\tau_2)$-stable,
\item all cocycle triples $(\dbar_1,\dbar_2,\phi)\in
\pi^{-1}(\dbar_1,\dbar_2,\Phi)$\ are
$(a_1,a_2,\tau_1,\tau_2)$-stable,
\item given any $(\dbar_1,\dbar_2,\phi)\in
\pi^{-1}(\dbar_1,\dbar_2,\Phi)$, every cocycle triple on
the
$\Omega^{p-1}_{1,2}$-orbit through
$(\dbar_1,\dbar_2,\phi)$\
is
$(a_1,a_2,\tau_1,\tau_2)$-stable.
\endroster
\endproclaim

\demo{Proof} Statements (2) and (3) are obviously
equivalent. We thus need only
prove that (2) or (3) is equivalent to (1).  To do so,
we need to compare the definitions of stability for
a cocycle triple and for a cohomology triple.  In both
cases,
the definition is given in terms of the values of
$\theta_{a_1,a_2,\tau_1,\tau_2}(\Cal E'_1,\Cal E'_2)$,
where $\E'_1$\ and $\E'_2$\ are the subbundles
in either a cocycle subtriple, $(\dbar'_1,\dbar'_2,\phi')$,
or
a cohomology subtriple, $(\dbar'_1,\dbar'_2,\Phi')$.
Notice that neither the $\phi'$\ nor the $\Phi'$\ affect
the
value of $\theta_{a_1,a_2,\tau_1,\tau_2}$ - their only
role is to determine \it on which pairs\rm\
$(\Cal E'_1,\Cal E'_2)$\ the function must be evaluated.
The proof thus consists essentially of a comparison of
the \it subobjects\rm\   of cocycle triples and of
cohomology triples.

Suppose first that $(\dbar_1,\dbar_2,\Phi)\in\Cal
H^{(p)}_{coh}$\
is
$(a_1,a_2,\tau_1,\tau_2)$-stable.  Let
$(\dbar_1,\dbar_2,\phi)$\ be any
cocycle triple in $\pi^{-1}(\dbar_1,\dbar_2,\Phi)$, and
let
$(\dbar'_1,\dbar'_2,\phi')$\ be a cocycle subtriple. Then
$\phi'$\ defines a cohomology
class, $\Phi'$, in $ H^1(X,Hom(\E_2,\E_1)$, and
$(\dbar'_1,\dbar'_2,\Phi')$\ is clearly a cohomology
subtriple of
$(\dbar_1,\dbar_2,\Phi)$.  Thus, by the stability of
$(\dbar_1,\dbar_2,\Phi)$,
$$\theta_{a_1,a_2,\tau_1,\tau_2}(\Cal E'_1,\Cal E'_2)<
0\
,$$
i.e. $(\dbar_1,\dbar_2,\phi)$\ is
$(a_1,a_2,\tau_1,\tau_2)$-stable.

Conversely, suppose that all cocycle triples,
$(\dbar_1,\dbar_2,\phi)$, in
$\pi^{-1}(\dbar_1,\dbar_2,\Phi)$\ are
$(a_1,a_2,\tau_1,\tau_2)$-stable.  Let
$(\dbar'_1,\dbar'_2,\Phi')$\ be any cohomology subtriple.
Fix any
representatives $\phi'$( resp. $\phi$) of $\Phi'$(
resp.$\Phi$).  Then
the condition $i_*(\Phi')=r_*(\Phi)$\ implies that
$i(\phi')= r(\phi)+(\dbar_1\circ\alpha'-
\alpha'\circ\dbar'_2)$,
for some $\alpha'\in\Omega^0(X,Hom(E'_2,E_1)$.  Let
$\alpha\in\Omega^0(X,Hom(E_2,E_1)$\ be any element
such that $\alpha\circ j=\alpha'$,
where $j:\E'_2\hookrightarrow\E_2$\ is the inclusion
map.  Then, since $j$\ is a holomorphic map, we get
$$\align
\dbar_1\circ\alpha'-\alpha'\circ\dbar'_2&=
\dbar_1\circ\alpha\circ j-\alpha\circ j\circ\dbar'_2\\
&=\dbar_1\circ\alpha\circ j-\alpha\circ \dbar_2\circ
j\\
&= r(\dbar_{1,2}(\alpha))
\endalign$$
That is, $i(\phi')= r(\phi+\dbar_{1,2}(\alpha))$, and hence
$(\dbar'_1,\dbar'_2,\phi')$\ is a cocycle subtriple
of $(\dbar_1,\dbar_2,\phi+\dbar_{1,2}(\alpha))$.
Since $\pi(\dbar_1,\dbar_2,\phi+\dbar_{1,2}(\alpha))
=(\dbar_1,\dbar_2,\Phi)$, it now follows from the
stability
of this cocycle
triple that $\theta_{a_1,a_2,\tau_1,\tau_2}(\Cal E'_1,\Cal E'_2)< 0$,
i.e. $(\dbar_1,\dbar_2,\Phi)$\ is
$(a_1,a_2,\tau_1,\tau_2)$-stable.
\hfill \qed
\enddemo
\bigskip

\subheading{\hfil\S 2.4 Metric equations}

In this section we describe the  metric equations
corresponding to the above definitions of stability.
Recall that for a triple $\tri$\ with $\Phi\in
H^0(X,Hom(E_2,E_1)$, there is a Hitchin-Kobayashi
correspondence between stability (as defined in [BGP])
and
metrics satisfying the coupled vortex equations. As
equations for metrics $H_1$\ and $H_2$\ on $E_1$\ and
$E_2$,
these are

$$i\Lambda F_{H_1}+\Phi\Phi^*=\tau_1\bold I\
,\tag 2.4.1a$$
$$i\Lambda F_{H_2}-\Phi^*\Phi=\tau_2\bold I\
.\tag 2.4.1b$$
where $\Phi^*$\ denotes the adjoint with respect to the
metrics $H_1$\ and $H_2$. To obtain the analogous
equations
corresponding to $(a_1,a_2,\tau_1,\tau_2)$-stability of
a
p-cohomology triple, we need the following operations on
form-valued sections of bundles over \kler\ manifolds
(cf
[W]).
$$\wedge:\Omega^{p,q}(X,E)\times
\Omega^{k,l}(X,E^*)\-->
\Omega^{p+k,q+l}(X,\Bbb C)\ ,\tag 2.4.2$$
$$\circ : \Omega^{p,q}(X,Hom(E_1,E_2))\times
\Omega^{k,l}(X,Hom(E_2,E_1))\-->
\Omega^{p+k,q+l}(X,Hom(E_1,E_1))\ ,\tag 2.4.3$$
$$\*E:\Omega^{p,q}(X,E)\-->
\Omega^{n-p,n-q}(X,E^*)\ .\tag 2.4.4$$
These are defined such that for $\varphi_i\in\Omega^{p,q}(X,E)$,

$$\varphi_1\wedge\*E\varphi_2=
(\varphi_1,\varphi_2)\frac{\omega^n}{n!}\tag 2.4.5$$

where $\omega$\ is the \kler\ form, and $(\phi,\psi)$\
is
the inner product coming from the metric on $E$\ and the
metric on forms of type (p,q). Also, for
$\phi_i\in\Omega^{p,q}(X,Hom(E_1,E_2))$, we have

$$\phi_1\wedge\*E \phi_2=Tr(\phi_1\circ\*E \phi_2)\tag 2.4.6$$

\proclaim{Definition 2.11}Given a p-cohomology triple $\tri$, and real
parameters
$(a_1,a_2,\tau_1,\tau_2)$\ we define the
following
equations for metrics on $E_1$\ and $E_2$\ and a
representative $\phi\in\Omega^{0,p}(X,Hom(E_2,E_1))$\
of the
cohomology class $\Phi$:

$$i\Lambda a_1F_{H_1}+\Lambda^n(\phi\circ\*E\phi)=
\tau_1\bold I\ ,\tag 2.4.7a$$
$$i\Lambda a_2 F_{H_2}-(-1)^p
\Lambda^n(\*E\phi\circ\phi)=
\tau_2\bold I\ ,\tag 2.4.7b$$
$$\dbar_{1,2}^*(\phi)=0\ .\tag 2.4.7c$$
\endproclaim

\noindent\bf Remarks 2.12.\rm\

\bf 2.12.1\rm\ The sign of the terms involving $\phi$\ are chosen
such that $Tr(\Lambda^n(\phi\circ\*E\phi))$\ and
$Tr((-1)^p\Lambda^n(\*E\phi\circ\phi))$\ are positive.
This will be important in section 2.6.

\bf 2.12.2\rm\ The coefficients $a_1$\ and
$a_2$\ will be assumed non-negative, and  the  parameters
$(a_1,a_2,\tau_1,\tau_2)$\ must satisfy
the constraint
$$a_1d_1+a_2d_2-\tau_1r_1-\tau_2r_2=0\ .\tag 2.4.8$$

\bf 2.12.3\rm\ The coupled vortex equations
given in [BGP] correspond to the case $p=0$\ and
$a_1=a_2=1$.  There is however no good reason to single
out
these special values for $a_1,a_2$. This is most clearly
seen in the symplectic interpretation of the equations,
and
will be discussed in the next section. We remark in passing that there is
no need to add scale factors to the terms involving
$\phi$ since these can be absorbed in $\phi$, or by a rescaling of the
metrics.

\subheading {\hfil \S 2.5 Moment maps}

If we fix metrics $K_1$\ and $K_2$\ on $E_1$\ and
$E_2$, we
can
reduce the gauge groups to the real unitary groups
$\G_1$\
and $\G_2$.
In addition,
$\Cal C_i,\ i=1,2$\ and $\Omega^{0,p}(X,Hom(E_2,E_1)$\
acquire
symplectic structures in the usual way. We denote these
by
$\omega_1$,
$\omega_2$, and $\omega_{(0,p)}$\ respectively.  A
symplectic form on
$$\chi^{(p)}
=\C_1\times\C_2\times\Omega^{0,p}(Hom(E_2,E_1)\
,$$
can be produced
by taking the sum
$\omega_1+\omega_2+\omega_{(0,p)}$.  This is,
however, merely one possibility; given any real
positive numbers $a_1$\ and $a_2$, we can form
a symplectic structures on $\chi^{(1)}$\ by defining
$$\omega_{a_1,a_2}=a_1\omega_1+a_2\omega_2
+\omega_{(0,p)}\ .\tag 2.5.1$$

\proclaim{Lemma 2.13}
The group $\G_1\times\G_2$\ acts symplectically on
$(\chi,\omega_{a_1,a_2})$, and has a moment map
$$\Psi_{a_1,a_2}:\chi\-->\lieG_1\times\lieG_2\ $$
given by
$$\Psi_{a_1,a_2}(\dbar_1,\dbar_2,\phi)=
( a_1\Lambda F_{K_1}-i\Lambda^n(\phi\circ\*E\phi),
 a_2\Lambda F_{K_2}-i(-1)^p
\Lambda^n(\*E\phi\circ\phi)\ .\tag 2.5.2$$
\endproclaim

\demo{proof} Exactly the same as for the $p=0$\ case.
The
sign factor $(-1)^p$\ comes from interchanging the order
in a wedge product forms of type (0,p) and (n, n-p).
\enddemo

If we define $\Cal
H^{eqtn}_{a_1,a_2,\tau_1,\tau_2}\subset\Cal
C_1\times\Cal
C_2\times\Omega^{0,p}(X,Hom(E_2,E_1))$\ to be the set
of
triples $(\dbar_1,\dbar_2,\phi)$\ on which solutions can
be
found to the coupled equations (a) and (b), then we have

\proclaim{Proposition 2.14}There is a bijective
correspondence
$$\Cal
H^{eqtn}_{a_1,a_2,\tau_1,\tau_2}/\Gc^{(1)}\times\Gc^{(2
)}
\longleftrightarrow (\Psi^{-1}_{a_1,a_2}(-i\tau_1,-
i\tau_2)
\cap\Cal H)/\G_1\times\G_2\ .$$
\endproclaim

Unfortunately, there does not seem to be a way to realize
the harmonicity condition as a moment map condition.

\subheading{\hfil \S 2.6 Hitchin-Kobayashi
correspondence}

In this section we show how our
stability conditions follow as a consequence from the
existence of solutions to the appropriate metric
equations.

\proclaim{Lemma 2.15}
Let $(\E_1,\E_2,\phi)$\ be a p-cocycle triple (so
$\phi\in\Omega^{0,p}(Hom(\E_2,\E_1))$
and $\dbar(\phi)=0$). Let $(a_1,a_2,\tau_1,\tau_2)$\ be
any
set of real numbers with $a_i\ge 0$\ for $i=1,2$, and
such
that
$$a_1d_1+a_2d_2-\tau_1r_1-\tau_2r_2=0\ .$$
Suppose there are bundle metrics $H_1$\ and $H_2$\
which
satisfy
the coupled equations (2.4.7a,b), i.e.
$$i\Lambda a_1F_{H_1}+\Lambda^n(\phi\circ\*E \phi)=
\tau_1\bold I\ ,$$
$$i\Lambda a_2F_{H_2}-(-1)^p
\Lambda^n(\*E\phi\circ\phi)
=\tau_2\bold I\ .$$

Let $(\E'_1,\E'_2,\phi')$
be a locally free subtriple, i.e. suppose that
$i_1:\E'_1\hookrightarrow\E_1$\ is
a subbundle of $\E_1$, $\E'_2\hookrightarrow\E_2$\ is
a subbundle of $\E_2$, and
$\phi'\in\Omega^{0,1}(Hom(\E'_2,\E'_1))$\ satisfies the
condition $i_1\circ\phi'=\phi\circ i_2$.

Then
$$\Theta_{a_1,a_2,\tau_1,\tau_2}(\E'_1,\E'_2)
\le 0\ ,$$
with equality if and only if $(\E_1,\E_2,\phi)$\ splits
with $(\E'_1,\E'_2,\phi')$\ as a direct summand.
\endproclaim

\demo{Proof}
Using the
metrics $H_1$\ and $H_2$\ we
can make orthogonal decompositions $\E_1=\E'_1\oplus
(\E_1/\E'_1)$\ and
$\E_2=\E'_2\oplus (\E_2/\E'_2)$. With respect to these
decompositions we can write
$\phi$\ as
$$\phi=\pmatrix
\phi'' & \phi^{'\perp} \\
\phi^{\perp '}& \phi^{\perp\perp}
\endpmatrix\ .\tag 2.6.1$$
However,  the definition of a subtriple requires that
$r(\phi)=i(\phi')$, where
$r$\ and $\phi$\ are the maps in
$$\CD Hom(\E_2,\E_1)@>r>> Hom(\E_2',\E_1) @<i<<
Hom(\E'_2,\E'_1)\endCD\ .$$
It follows that in (2.6.1) we have $\phi''=\phi'$\ and
$\phi^{\perp '}=0$, i.e.
$$\phi=\pmatrix
\phi' & \phi^{'\perp} \\
0 & \phi^{\perp\perp}
\endpmatrix\ .\tag 2.6.2$$
The conclusion now follows precisely as in the case of
ordinary triples.  More specifically, after writing the
curvature terms with respect to the above orthogonal
decompositions of the bundles, the equations (2.4.7a,b) yield the following:

$$ia_1\Lambda F'_{H_1}+ia_1\Lambda\Pi_1+
 \Lambda^n(\phi''\circ\*E\phi'')+
\Lambda^n(\phi^{'\perp}\wedge\*E \phi^{'\perp})=
\tau_1\bold I'_1\ , \tag 2.6.3a$$
$$ia_2\Lambda F'_{H_2}+ia_2\Lambda\Pi_2-
(-1)^p\Lambda^n ( \phi''\wedge\*E\phi'')=
\tau_2\bold I'_2\ , \tag 2.6.3b$$

We can take the trace of these equations, and use the
fact that for any section
$\phi\in\Omega^{0,p}(Hom(\E_2,\E_1))$\ we have
$$Tr((-1)^p\Lambda^n\*E\phi\wedge\phi)=
Tr(\Lambda^n\phi\wedge \*E\phi)=|\psi|^2\ .\tag 2.6.4$$
This gives

$$a_1d'_1+a_2d'_2
+a_1Tr(i\Lambda\Pi_1)+a_2Tr(i\Lambda\Pi_2)
+|\phi^{'\perp}|^2=\tau_1 r'_1+\tau_2r'_2\ .\tag 2.6.5$$
The conclusion follows directly from this, since both
$Tr(i\Lambda\Pi_1)$\ and $Tr(i\Lambda\Pi_2)$\ are
non-
negative. \hfill\qed
\enddemo

We now consider p-cohomology triples over Riemann
surfaces

\proclaim{Theorem 2.16}

Let $(\E_1,\E_2,\Phi)$\ be a p-cohomology triple over a
Riemann surface $X$ (so
$\Phi\in H^p(Hom(\E_2,\E_1))$).
Let $(a_1,a_2,\tau_1,\tau_2)$\ be any set of real
numbers with $a_i\ge 0$\ and
satisfying the constraint
$$a_1d_1+a_2d_2-\tau_1r_1-\tau_2r_2=0\ .$$
Suppose there is a representative
$\phi\in\Omega^{0,p}(Hom(\E_2,\E_1))$\ for $\Phi$, and
bundle metrics $H_1$\ and $H_2$, which satisfy
the coupled equations (2.4.7a-c).  Then
$(\E_1,\E_2,\Phi)$\
is
$(a_1,a_2,\tau_1,\tau_2)$-stable.
\endproclaim

\demo{Proof}
Since $X$ is a Riemann surface, all subtriples are locally
free.  Let $(\E'_1,\E'_2,\Phi')$\ be any such
subtriple. To prove the Proposition we need to
show that
$$\theta_{a_1,a_2,\tau_1,\tau_2}(\Cal
E'_1,\Cal E'_2)< 0\ .$$

Notice that by equation (2.4.7c), $\phi$\ is the harmonic
representative of $\Phi$\ with respect to the metrics
$H_1$\
and $H_2$. Using the induced metrics
on $\E'_1$\ and $\E'_2$, take the harmonic
representative,
$\phi'$, of
$\Phi'$. We claim that $(\E'_1,\E'_2,\phi')$\ is a
subtriple of the 1-cocycle triple $(\E_1,\E_2,\phi)$, i.e.
we claim that $r(\phi)=i(\phi')$.  This will prove the proposition, since then
by Lemma 2.12, we have
$\theta_{a_1,a_2,\tau_1,\tau_2}
(\Cal E'_1,\Cal E'_2)< 0$.

We now prove our claim. From the very definition of the
maps induced by $r$\ and $i$\ in cohomology, we get that
$r(\phi)=i(\phi')+\dbar(\alpha)$,
where $\alpha\in\Omega^0(Hom(\E'_2,\E_1))$\ and
$\dbar$\
denotes the operator
induced by $\dbar'_2$\ and $\dbar_1$. With respect to
the orthogonal decompositions $\E_1=\E'_1\oplus (\E_1/\E'_1)$\ and
$\E_2=\E'_2\oplus (\E_2/\E'_2)$\ we can thus write

$$\phi=\pmatrix
\phi'+\beta'' & \phi^{'\perp} \\
\beta^{\perp '}& \phi^{\perp\perp}
\endpmatrix\ ,$$
where $\dbar(\alpha)=\beta''+\beta^{\perp '}$. But
harmonic
representative
are norm minimizing. Thus $||\phi'+\beta''||^2\ge
||\phi'||^2$, and therefore
$$
||\phi||^2\ge \Vmatrix
\phi' & \phi^{'\perp} \\
0 & \phi^{\perp\perp}
\endVmatrix ^2\ .$$
Thus we get $||\phi||^2\ge ||\phi-\dbar(\alpha)||^2$, which
is a contradiction unless $\dbar(\alpha)=0$.\qed
\enddemo

\heading\S 3 Extensions\endheading

The case of 1-cohomology triples deserves special attention because of the
fact
that a 1-cohomology class
$\Phi\in H^1(Hom(\E_2,\E_1))$\ can be interpreted as an
extension class for extensions of $\E_2$\ by $\E_1$.
This can be exploited to study moduli space questions
for the set of all such extensions, i.e. for the set of all
short exact sequences
$$\extn\ ,\tag e$$
where $\E_1$\ and $\E_2$\ have fixed underlying smooth
bundles (denoted by $E_1$\ and $E_2$\ respectively).

Such extensions can also be considered from the point of
view of the bundle $\Cal E$.  This leads to a metric
problem and definition of stability that appear somewhat
different to the ones considered in the previous section.
In this section we discuss such an approach. In the next
section we indicate the relationship between the two
approaches.

Let us begin therefore with a compact \kler\ manifold
$X$,
and a holomorphic bundle $\Cal E\-->X$\ given as an
extensions of bundles as in (e).

\subheading{\hfil\S 3.1 stability}

To formulate the stability condition, we consider extensions as objects in the
category of short exact sequences of coherent sheaves of the form

$$\fextn\ . \tag f $$

A morphism of two extensions is defined by the commutative diagram
$$\CD
 0 @>>>   \cF_1   @>>> \cF  @>>>  \cF_2 @>>>  0\\
   @.      @A{f_1}AA @A{F}AA    @AA{f_2}A    @.    \\
 0 @>>>    \cF_1' @>>>   \cF'  @>>>    \cF_2' @>>>    0.
\endCD\tag 3.1.1$$
\noindent {\bf Remark 3.1}.\
This category is abelian. In particular,
$$
0 \lra  \ker f_1  \lra  \ker f    \lra   \ker f_2    \lra     0,
$$
and
$$
 0 \lra  \im  f_1  \lra  \im f    \lra   \im f_2    \lra     0,
$$
are the kernel and image of the map $F$\ in (3.1.1).

A subobject of (f) consists then of an extension
$$
 0 \lra  \cF_1'  \lra  \cF'    \lra   \cF_2'    \lra     0,
\tag f'$$
and injective maps $i_1$, $i_2$, $i$ such that the following diagram
commutes
$$\CD
 0 @>>>   \cF_1   @>>> \cF  @>>>  \cF_2 @>>>  0\\
   @.      @A{i_1}AA @A{i}AA    @AA{i_2}A    @.    \\
 0 @>>>    \cF_1' @>>>   \cF'  @>>>    \cF_2' @>>>    0.
\endCD\tag 3.1.2$$
The extension (f') will be called a \it subextension\rm\ of (f).

\proclaim{Lemma 3.2}
Let us consider the extension (f).
Any subsheaf $i:\cF'\hookrightarrow\cF$ defines a subextension of (f).
\endproclaim
\demo{Proof}
Let $g:\cF'\-->\cF_2$ be the map obtained by composing $i$ with the surjection
$\cF\-->\cF_2$. Then
$$
0 \lra  \ker g  \lra  \cF'    \lra   \im g    \lra     0
$$
is the desired subextension.\qed
\enddemo

There is hence a one-to-one correspondence between subsheaves of $\cF$ and
subextensions
of (f).

\proclaim{Definition 3.3}Let $e$ be  the extension
$$
\extn
$$
and  $e'$ the subextension
 $$
\sextn.
$$
For $\alpha\in \Bbb R$ we  define the $\alpha$-slope of $e'$ as
$$
\mu_\alpha(e')=\mu(\cE')+\alpha \; \frac{rank\cE_2'}{rank \cE'}.\tag 3.1.3 $$
\endproclaim

\proclaim{Definition 3.4}
The extension $e$ is said to be $\alpha$-stable (resp. semistable) if and
 only if
for every subextension $e'\subset e$ (resp. $e'\subseteq e$)
$$
\mu_\alpha(e')<\mu_\alpha (e)\;\;\; (\text{resp.}\;\; \leq).\tag 3.1.4
$$
\endproclaim

\noindent\bf Remark.\rm\ If $\alpha=0$, then \asty\ is equivalent to ordinary
stability.

\proclaim{Proposition 3.5}
Let $e$ be \as, then
$$
\alpha> \mu(\cE_1)-\mu(\cE_2).\tag 3.1.5
$$
\endproclaim
\demo{Proof}
It suffices to apply the numerical stability condition to the trivial
subextension
$$
0\lra \cE_1\lra \cE_1\lra 0\lra 0.
$$
\qed
\enddemo

On the other hand we will see later (when dealing with the Donaldson
functional) that $\alpha\leq 0$, and hence in order for $e$ to be \as\ it is
necessary that $\alpha$  lies in the interval
$$
(\mu(E_1)-\mu(E_2), 0].\tag 3.1.6
$$
As usual this interval will be subdivided by some critical points, and the
stability
condition will just depend on the subinterval. Moreover one has the following.

\proclaim{Lemma 3.6} There is some $\epsilon>0$\ such that
for $\alpha$\ in the interval
$$(\mu(E_1)-\mu(E_2),\mu(E_1)-\mu(E_2)+\epsilon)\
,$$
the following is satisfied:

(1)\ If $e$ is  \as, then $\cE_1$ and $\cE_2$ are semistable.

(2)\ If $\cE_1$ and $\cE_2$ are stable then $e$ is \as.
\endproclaim
\demo{Proof} The proof is essentially identical to that of Proposition 2.9. In
fact, in view of the results of section 4.2, this result can be treated as a
special case of Proposition 2.9, corresponding to the case $p=1, a_1=a_2=1$.
We can also give a direct proof which depends on an examination of the
\astability\ condition for special sub-objects. In this case the subobjects
are
subextensions with either $\cE'_2=0$\ or $\cE'_1=\cE_1$.\qed
\enddemo

\subheading{\hfil\S 3.2 metric equations}

Given an extension
$$\extn\ ,\tag e$$
the natural metric problem is to look for a metric $H$ on $\cE$ satisfying the
equation
$$
i\Lambda F_H=
\pmatrix
\tau_1 I_1& 0\\
0&\tau_2 I_2
\endpmatrix\ .\tag 3.2.1$$
Here $\tau_1$ and $\tau_2$ are real numbers an $I_1$ and $I_2$ are the
identity
endomorphisms in $E_1$ and $E_2$ respectively. We can make sense of the right
hand side since the metric $H$ on $\cE$
gives  a $C^\infty$ splitting of ($E$), i.e. an identification
of the smooth underlying bundle to $\cE$ with $E_1\oplus E_2$.

\noindent\bf Remark.\rm\  If $\tau_1=\tau_2=\lambda $, equation (3.2.1)
reduces
to the Hermitian-Einstein equation.

\proclaim{Proposition 3.7}
If $H$ satisfies (3.2.1), then the parameters $\tau_1$ and $\tau_2$ are
related
by
$$
r_1\tau_1+r_2\tau_2=d_1+d_2,
\tag 3.2.2$$
where $r_1=rank\cE_1$, $r_2=rank \cE_2$, $d_1=\deg\cE_1$ and $d_2=\deg\cE_2$.
\endproclaim
\demo{Proof}
This is easily proved by taking the trace in both sides of (3.2.1) and
integrating.
\enddemo

\subheading{\hfil \S 3.3 Hitchin-Kobayashi correspondence}

We first prove that \astability\ is a necessary condition for existence of
solutions to the equation (3.2.1).

\proclaim{Proposition 3.8}
Let $e$ be the extension of vector bundles
$$\extn,\tag e$$
and let $\tau_1$ and $\tau_2$ satisfy (3.2.2). Set $\alpha=\tau_1-\tau_2$.
If $\cE$ is indecomposable and admits a metric $H$ satisfying the metric
equation
for extensions (e), then $e$ is \as.
\endproclaim

\demo{Proof}
We need to show that $\mu_\alpha(e')<\mu_\alpha (e)$\ for every subextension
$e'$
$$\sextn\ .$$
The proof is a minor modification of the analogous result for the ordinary
Hitchin-Kobayashi correspondence.  Consider first the locally free
subextensions, i.e. the $e'$\ in which $\cE_1'$\ and $\cE_2'$\ are locally
free. Denote the underlying smooth bundle for $\cE'$\ by $E'$, and let
$E^{\perp}$\ be its orthogonal complement with respect to $H$. Then with
respect to the smooth orthogonal splitting $E=E'\oplus E^{\perp}$, we get the
block diagonal decomposition
$$\sqrt{-1}\Lambda F_{H}=
\pmatrix
 \sqrt{-1}\Lambda F'+\Pi' & \ast\\
 \ast& \sqrt{-1}\Lambda F^{\perp}-\Pi^{\perp}
\endpmatrix\tag 3.3.1$$
 where $\Lambda F'$\ and $\Lambda F^{\perp}$\ are the induced metric
connections on $\cE'$ and $\cE^{\perp}$\ respectively, and
$\Pi'$,$\Pi^{\perp}$\ are positive definite endomorphisms coming
from the second fundamental form for the inclusion of
 $E'$\ in $E$. With respect to this splitting of $E$, the endomorphism on the
right hand side of the metric equation is no longer diagonal, but has the form
$$\pmatrix
 T' & \ast\\
 \ast& T^{\perp}
\endpmatrix
=A\pmatrix
\tau_1 I_1& 0\\
 0&\tau_2 I_2
\endpmatrix A^{-1}\ ,\tag 3.3.2$$

where the matrix $A$\ gives the transformation from the frame $E_1\oplus E_2$\
to $E'\oplus E^{\perp}$. If we make the further orthogonal decompositions of
$E_1$\ and $E_2$\ into components in $E'$\ and $E^{\perp}$, then
$$E_1\oplus E_2=E_1'\oplus E_1^{\perp}\oplus
 E_2'\oplus E_2^{\perp}\ ,$$
and
$$E'\oplus E^{\perp}=E_1'\oplus
 E_2'\oplus E_1^{\perp}\oplus E_2^{\perp}\ .$$
With respect to these frames, the transformation $A$\ is represented by
 $$\pmatrix
 I_1' & 0 & 0 & 0 \\
 0 & 0 & I'_2 & 0 \\
 0 & I^{\perp}_1 & 0 & 0 \\
 0 & 0 & 0 & I^{\perp}_2 \\
\endpmatrix\tag 3.3.3$$
In fact, all we need is the trace of $T'$. It follows by a straightforward
linear algebra computation that
$$Tr(T')=r_1'\tau_1+r_2'\tau_2\ ,$$
where $r_1'=$rank$\cE_1'$\ and $r_2'=$rank$\cE_2'$.
We apply this to the condition
$$\sqrt{-1}\Lambda F'+\Pi'=T'\ ,$$
which can be extracted from the full metric equation. After taking the trace
and integrating over $X$, we thus get
$$\int_X{Tr(\sqrt{-1}\Lambda F')}+\int_X{Tr(\Pi')}=r_1'\tau_1+r_2'\tau_2\
.\tag
3.3.4$$
Using the Chern-Weil formula for deg$(\cE')$, and the positivity of $\Pi'$, we
obtain

$$deg(\cE')\le r_1'\tau_1+r_2'\tau_2\ ,\tag 3.3.5$$

with equality if and only if $\Pi'=0$, i.e. if and only if $\cE$\ splits. If
$\alpha=\tau_1-\tau_2$, and $\mu_{\alpha}(\cE')$\ is as in Definition 3.3,
then
(3.3.5) is equivalent to  $\mu_\alpha(e')<\mu_\alpha (e)$.
This proves the result for locally free subobjects.
If $e'$\ is not locally free, then there is a subvariety
$\Sigma\subset X$\ of codimension at most two, such that
$\cE'|_{X-\Sigma}$\ is locally free. We can thus apply the above arguments
over
$X-\Sigma$. This is good enough, because of the size of the codimension of
$\Sigma$.\qed
\enddemo

We now prove that \astability\ is a sufficient condition for existence of
special metrics as defined the deformation of the Hermitian-einstein equations
given in (3.2.1). That is, we prove

\proclaim{Theorem 3.9}
Suppose that $\alpha<0$\ and
$$\extn$$
 is an $\alpha$-stable extension. Let $\tau_1$\ and $\tau_2$\ be such that
$\alpha=\tau_1-\tau_2$\ and $deg(E)=r_1\tau_1+r_2\tau_2$.  Then there is a
metric $H$\ on $\cE$\ satisfying the equation (3.2.1), i.e
$$i\Lambda F_H = \pmatrix
                 \tau_1& 0\\
                  0 &\tau_2
\endpmatrix$$
\endproclaim

The proof is an adaptation of the methods used in [Do] (also [S] and [U-Y]) in
proving the Hitchin-Kobayashi correspondence for ordinary stable bundles.  As
shown by Donaldson, the \he\ equation is the equation satisfied by the
critical
points of a certain functional defined on the space of Hermitian metrics on
$\cE$. We shall modify this functional to show that our equations appear in
the
same way.

Just as in the case of the Hermitian-Einstein equation, we can separate out
the
Trace and Trace-free parts of the equation. We can fix the determinant of the
metric on $E$ to satisfy the trace part,
$$
i\Lambda Tr(F_H)=r_1\tau_1+r_2\tau_2.\tag 3.3.6
$$
The problem then becomes one of finding a new metric with this same
determinant, and which satisfies
$$i\Lambda F^0_H = \pmatrix
                 \frac{r_2}{r}\alpha& 0\\
                  0 &-\frac{r_1}{r}\alpha
\endpmatrix\tag 3.3.7$$
where $\alpha=\tau_1-\tau_2$.


Recall Donaldson's original functional to prove existence of solutions of the
Hermitian-Einstein\ equation:
Let $\cE$ be a holomorphic vector bundle over a compact \kler\ manifold
$(X,\omega)$. Donaldson defined a functional $M(-,-)$
on pairs of Hermitian metrics on $\cE$ using Bott--Chern secondary classes.
Namely
$$ M(H,K)=\int_X(R_2(H,K)-2\lambda R_1(H,K)\omega)\wedge\omega^{n-1},\tag
3.3.8$$
where
$$\align
R_1(H,K)&=\log\det(K^{-1}H)=\tr(\log K^{-1}H)\tag 3.3.9a\\
i\dbar\partial R_2(H,K)&=(-\tr(F_H^2))-(-\tr(F_K^2)),\tag 3.3.9b\\
\lambda &= \frac{\deg \cE}{rank \cE}\ .\tag 3.3.9c
\endalign $$

Now fix a smooth background metric $K$, with determinant satisfying (3.3.6).
Let
$$S(K)=\{s\in\Omega^0(X,EndE) | s^{*_K}=s\ ,\ Tr(s)=0\}
\ .\tag 3.3.10$$
Then any other metric with the same determinant as $K$\ can be described by
$Ke^s$, with $s\in S(K)$.  Fix an integer $p>2n$, and define
$$\Met^p_2=\{H=Ke^s \ s\in L^p_2(S(K))\}\ .\tag 3.3.11$$

Let $M:\Met(\cE)\-->\Bbb R$ be given by
$M(H)=M(K,H)$. The important property of $M$ is that $H$ is a critical
point if and only if $H$ satisfies the trace free part of the \he\ equations ,
i.e.
$$
i\Lambda F^0_H=0.
$$

Consider now the extension

$$\extn\ .\tag e$$
Given a background metric $K$ on $\cE$ we can (smoothly) identify $\cE_2$\
with
the orthogonal complement of $\cE_1$\ in $\cE$, and in this way get metrics
$K_1$ and $K_2$ on $\cE_1$ and $\cE_2$ respectively.
Any other metric H can similarly be split into $H_1$\ and $H_2$\ (by using an
$H$-orthogonal splitting of $\cE$). Denote
$$
M_D(H,K)=\int_X R_2(H,K)\wedge\omega^{n-1}.
\tag 3.3.12$$
Let $\tau_1$ and $\tau_2$ be real parameters. We shall consider the functional
$$
M_{\tau_1,\tau_2}(H,K)=M_D(H,K)-2\int_X(\tau_1 R_1(H_1,K_1) + \tau_2
R_1(H_2,K_2))
\wedge\omega^n.\tag 3.3.13
$$
\bf Remark.\rm\
If $\tau_1=\tau_2=\lambda$, then $
M_{\tau_1,\tau_2}(H,K)=M(H,K)$,
as can be easily seen from the following simple fact (see [Do, Prop. 7. p.
10]).

\proclaim{Lemma 3.10}
Let $H$ and $K$ be Hermitian metrics on $\cE$ and Let $H_1$, $K_1$ and
$H_2$, $K_2$ the corresponding metrics induced on $\cE_1$ and $\cE_2$
respectively, then
$$
R_1(H,K) =R_1(H_1,K_1) + R_1(H_2,K_2).
$$
\endproclaim

Notice that $R_1(H,K)=0$\ if the metrics have fixed determinant.  We can thus
simplify our definition to
$$
M_{\tau_1,\tau_2}(H,K)=M_D(H,K)-2(\tau_1-\tau_2)\int_X R_1(H_1,K_1)
\wedge\omega^n\ .\tag 3.3.14
$$

Let us fix $K$ and define
$$
M_{\tau_1,\tau_2}(H)=M_{\tau_1,\tau_2}(H,K).
\tag 3.3.15$$

Define $m^0:\Cal Met\longrightarrow\Omega^0(X,EndE)$\ by
$$m^0(H)=\Lambda\curvdh^0
+\sqrt{-1}T^0_H\ ,\tag 3.3.16$$
where, with respect to the orthogonal splitting $E=E_1\oplus E_2$\ determined
by $H$,
$$\align
T^0_H&=\pmatrix
                 \tau_1 I_1& 0\\
                  0&\tau_2 I_2
                            \endpmatrix
-Tr\pmatrix
                 \tau_1 I_1& 0\\
                  0&\tau_2 I_2
                            \endpmatrix\\
&= \pmatrix
                 \frac{r_2}{r}\alpha& 0\\
                  0 &-\frac{r_1}{r}\alpha
\endpmatrix
\endalign$$

The crucial properties of $M_{\tau_1,\tau_2}$\ are described in the next
proposition.
\proclaim{Proposition 3.11}

(1) Given any three metrics $H,K,J$, we have
$$
M_{\tau_1,\tau_2}(H,K) +M_{\tau_1,\tau_2}(K,J)=M_{\tau_1,\tau_2}(H,J),
$$

(2)\ If $H(t)=He^{ts}$\ with $s\in S(H)$, then
$$
\frac{d}{dt} M_{\tau_1,\tau_2}(H(t))=
2i \int_X \tr\left(s
m^0(H(t))\right).$$

(3)\ If $s\in S(H)$\ is given by $s=
\pmatrix
                 s_1& u\\
                  u^*&s_2
\endpmatrix$
with respect to the orthogonal splitting $E=E_1\oplus E_2$\ determined by $H$,
then

$$\align
\frac{d^2}{dt^2} M_{\tau_1,\tau_2}(H(t))|_{t=0}&=
2i \int_X \tr\left(s\frac{d}{dt}m^0(H(t))|_{t=0}\right)\\
&=\parallel D'_H(s)\parallel^2-\alpha\parallel u\parallel^2
\endalign $$
\endproclaim

\demo{Proof}

(1)  This follows immediately from the properties of the Bott-Chern classes.

(2)  Chose a frame for $E$\ such that $H$ can be written
as
$$
\pmatrix H_1 & 0\\
         0  &  H_2
\endpmatrix\ ,$$
that is a frame in which $\cE_1$ and $\cE_2$ are $H$-orthogonal.  In terms of
this frame we can write
$$
s=\pmatrix s_1 & u\\
u^\ast  &  s_2
\endpmatrix\ ,$$
where $s_1\in S(H_1)$, $s_2\in S(H_2)$ and $u\in \Hom(E_2,E_1)$. We have to
show that
$$\align
\frac{d}{dt} M_{\tau_1,\tau_2}(H(t))|_{t=0}&=
2i \int_X \tr\left(s
i\Lambda F^0_H-
\pmatrix \frac{r_2}{r}\alpha I_1& 0\\
0& -\frac{r_1}{r}\alpha I_2\endpmatrix\right)\\
&=2i\int_X\tr\left(s\Lambda F^0_{H(t)}\right)
-2i\alpha\int_X\tr(s_1)\ .
\endalign $$
{}From [Do] we know that $
  \frac{d}{dt}M_D(H(t))=2i\int_X\tr\left(s\Lambda F^0_{H(t)}\right)$, so it
remains to compute
$\frac{d}{dt}R_1(H_1(t),H_1)|_{t=0}$.  If we write
$$H(t)=He^{ts}=H\pmatrix h_1(t) & *\\
*  &  *
\endpmatrix\ ,$$
with respect to the H-orthogonal frame, then
$$R_1(H_1(t),H_1)=\text{log det} h_1(t)=\text{log det}
(1+ts_1+\frac{t^2}{2}(s_1^2+uu^*)+O(t^3))\ .\tag 3.3.17$$
A straightforward computation yields the result
$$\frac{d}{dt}R_1(H_1(t),H_1)|_{t=0}= Tr(s_1)\ .$$

(3) It follows from  (400) that $$\frac{d^2}{dt^2}R_1(H_1(t),H_1)|_{t=0}
=Tr(uu^*)=|u|^2\ .$$
The result now follows from this, plus the fact that
$$\frac{d}{dt}M_D(H(t))=\parallel D'_H(s)\parallel^2\ .$$
\qed
\enddemo

Notice that as a consequence of (1) and (3) in Proposition 3.11 we get

\proclaim{Proposition 3.12} Suppose that $\alpha<0$\ and (e) is an
$\alpha$-stable extension.  Then
\roster
\item $$\frac{d^2}{dt^2} M_{\tau_1,\tau_2}(H(t))>0$$

\item $Ker(L)=0$, where $L$\ is the operator on $L^p_2(S(H)$\ defined by
$L(s)=\frac{d}{dt}m^0(H(t))|_{t=0}$.
\endroster
\endproclaim

\demo{Proof} Both of these statements follow from the fact that
if $s$\ is as in (3) and $L(s)=0$, then $\dbar_1(s_1)=\dbar_2(s_2)=u=0$. The
eiegenspaces of $s$\ thus split the extension (e) into a direct sum of
extensions. This violates the stability criterion, since the $\alpha$-slope
inequality cannot be satisfied by both
summands.\qed
\enddemo

The functional $M_{\tau_1,\tau_2}$\ thus has the convexity features we
require.
 Furthermore,

\proclaim{Lemma 3.13}
Suppose that $\alpha<0$\ and let $H=Ke^s$\ with $s\in L^p_2(S(K))$. Let $s=
\pmatrix s_1& u\\
u^*&s_2 \endpmatrix$\ be the block decomposition of $s$\
with respect to the orthogonal splitting $E=E_1\oplus E_2$\ determined by $K$.
Let
$\Psi:\Bbb R\times\Bbb R\longrightarrow\Bbb R$\ be the smooth function as in
[B] (or [S]). Then

$$\align
M_{\tau_1,\tau_2}(H)=\ &\sqrt{-1}\int_X{Tr(s\Lambda F_K)}+
\int_x{(\Psi(s)\dbare s,\dbare s)_K}
-2\alpha R_1(H_1,K_1)\\
\ge\ &\sqrt{-1}\int_X{Tr(s\Lambda F_K)}+
\int_x{(\Psi(s)\dbare s,\dbare s)_K}-\alpha\int_x{Tr(s_1)}
\tag 3.3.18 \endalign $$

where the meaning of $\Psi(s)$\ is as in [B]\ or [S].
\endproclaim

\demo{Proof} The first line follows from the computations in [S] (or [Do]).
The
second uses the convexity properties of the function $R_1(H(t)_1,K_1)$, and
the
fact that its first derivative at $t=0$\ is given by $\int_x(Tr(s_1)$.
\enddemo

This is slightly weaker than the analogous result for the original Donaldson
functional, but is strong enough for our purposes.

The rest of the proof of Theorem 3.9 is precisely along the lines of the
analogous result in [S].  We give here a sketch of the main ideas. Fix a real
number $B$\ such that
$\parallel m^0(K)\parallel^p_{L^p}\le B$ (where
$\parallel m^0(K)\parallel^p_{L^p}=\int_{X}{|m^0(K)|^p_K
dvol}$).  Define
$$\Cal Met^p_2(B)=\{H\in\Cal Met^p_2 |\
\parallel m^0(H)\parallel^p_{L^p}\le B\ \}\ $$

We look for minima of $M_{\tau_1,\tau_2}(H)$\ on $\Cal Met^p_2(B)$.  As in the
case of the unmodified Donaldson functional, if the extension (e) is
$\alpha$-stable, then there are no extrema on the boundary of this constrained
space, and  the minima occur at solutions to the metric equation $m^0(H)=0$.

To show that minima do occur, we need

\proclaim{Proposition 3.14} Either (e) is \as\ or we can find positive
constants
$C_1$\ and $C_2$\ such that
$$sup|s|<C_1 M_{\tau_1,\tau_2}(Ke^s)+C_2\ $$
for all $Ke^s\in \Met^p_2(B)$,
\endproclaim

\demo{sketch of Proof}As in the case of the unmodified Donaldson functional,
one first shows that for metrics in the constrained set $\Met^p_2(B)$, the
$C^0$ estimate given above is equivalent to a $C^1$\ estimate of the same
type.
One then supposes that no such estimate holds. It follows that one may find a
sequence
$\{u_i\}\subset L^P_2(S(K)$\ such that $\parallel u_i\parallel_{L^1}=1$. This
has a weakly convergent subsequence in  $L^2_1(S(K)$, with non-trivial limit
denoted by
$u_{\infty}$. One then shows that the eigenvalues of $u_{\infty}$ are constant
almost everywhere.  This is done, as in [S], by making use of an estimate of
the form:

\proclaim {Proposition 3.15} Let $\Cal F:\Bbb R\times \Bbb R\-->\Bbb R$\ be
any
smooth positive function which satisfies $\Cal F(x,y)\le 1/(x-y)$\ whenever
$x>y$.  Then
$$\sqrt{-1}\int_X{Tr(u_{\infty}\Lambda F_K)}+
\int_x{(\Cal F(u_{\infty})\dbare u_{\infty},\dbare
u_{\infty})_K}-\alpha\int_x{Tr(u_{\infty,1})}\le 0\ ,\tag 3.3.18$$
where $u_{\infty}=\pmatrix u_{\infty,1}& *\\
*&*\endpmatrix$\ with respect to the splitting of $E$ determined by $K$.
\endproclaim

\demo{Proof} This follows from the analysis in [S], plus the estimate given in
Lemma 3.3.13.
\enddemo

Since $Tr(u_{\infty})=0$, there are at least two distinct eigenvalues. Let
$\lambda_1<\lambda_2,\dots,<\lambda_k$\ denote the distinct eigenvalues.
Setting $a_i=\lambda_{i+1}-\lambda_i$, one can thus define projections
$\pi_i\in L^2_1(S(K))$\ such that

$$u_{\infty}=\lambda_r\bold I-\sum_{i}^{k-1}{a_i\pi_i}
\tag 3.3.19$$

By an important result of Uhlenbeck and Yau (cf. [U-Y]), the $\pi_i$\ define a
filtration of $\cE$\ by reflexive subsheaves
$$\cE_1\subset\cE_2\subset\dots\subset\cE_k=\cE$$
Each subsheaf $\cE_j$\ determines a subextension
$$ 0\longrightarrow\cE_{1,j}\longrightarrow\cE_{j}
\longrightarrow\cE_{2,j}\longrightarrow 0\ .$$
Now define the numerical quantity

$$Q=\lambda_k(r\mu(\cE)-r_1\tau_1-r_2\tau_2)-
\sum_{i}^{k_1}{a_i(r_i\mu(\cE_i)-
r_{1,i}\tau_1-r_{2,i}\tau_2)}\ ,\tag 3.3.20$$
where $\mu(\cE_i)$\ is the slope of $\cE_j$, and
$r_{a,i}$\ is the rank of $\cE_{a,i}$.

Using Lemma 3.3.13 and the fact that $u_{\infty}=\lambda_r\bold
I-\sum_{i}^{k-1}{a_i\pi_i}$, one shows (by precisely the method in [S]) that
$Q\le 0$.
On the other hand, $\tau_1$\ and $\tau_2$\ are related by
$r\mu(\cE)-r_1\tau_1-r_2\tau_2=0$, and if $(e)$\ is \as, then
$$ r_i\mu(\cE_i)-r_{1,i}\tau_1-r_{2,i}\tau_2<0$$
for all $i=1,\dots, k-1$. Thus $Q>0$\ if (e) is \as.
This proves the proposition.\qed
\enddemo

To complete the proof of the Hitchin-Kobayashi correspondence, it remains to
show that the minimum of
$M_{\tau_1,\tau_2}$\ is a smooth solution to the equation
(3.2.1). This is done exactly as in [Do] or [S].

\subheading{\hfil \S 3.4 An example}

As an example, we can consider the case where $E_1=L_1$\ and
$E_2=L_2$\ are line bundles over a Riemann surface.  We
assume further
that $d_1<d_2$, where $d_i$\ denotes the degree of $L_1$.
Let us denote such extensions by
$$0\-->\Cal L_1\-->\cE\-->\Cal L_2\-->0\ .\tag l$$

Since line bundles are automatically stable, Lemma 3.6(1)
gives

\proclaim{Lemma 3.16}Let $L_1$\ and $L_2$\ be as above.  Then
there is some
$\epsilon>0$\ such that all extensions $(l)$\ as above are
\astable, for
any $\alpha$\ in the interval $$(d_1-d_2,d_1-d_2+\epsilon)\ .$$
\endproclaim

We can give a more detailed analysis.  The main reason for
this is that the
possibilities for sub-extensions are so restricted; they all
correspond to
rank one (i.e. line-) subbundles of $\cE$, are are of one of
two types. The
only possibilities are
$$0\-->\Cal L_1\-->\Cal L_1\-->0\-->0\ ,\tag 3.4.1$$
or
$$0\-->0\-->\Cal L\-->\Cal L\-->0\ .\tag 3.4.2$$
Computing the $\alpha$-slopes, we see that
$$\align
\mu_{\alpha}(0,\Cal L_1)=&d_1\ ,\tag 3.4.3\\
\mu_{\alpha}(\Cal L,0)=&d_{\Cal L}+\alpha\ .\tag 3.4.4
\endalign$$
{}From this we see that if $(l)$\ is \astable, then
$$0<d_1-d_2<\alpha<d_1+d_2-2d_{\Cal L}\ \tag 3.4.5$$
for all subbundles $\Cal L\ne \Cal L_1$.  Now define
$$div(\Cal E) = Max\{d_{\Cal L}\ |\ d_{\Cal L}\ \text{is the degree of a line
subbundle of $\Cal E$}\}\ .\tag 3.4.6$$

\proclaim{Lemma 3.17}If
$$0<d_1-d_2<d_1+d_2-2div(\Cal E)\ ,\tag 3.4.7$$
then $(l)$\ is \astable\ for any $\alpha$\ in the interval
$(d_1-d_2\ ,\ d_1+d_2-2div(\Cal E))$.
\endproclaim
\demo{Proof} Given any $\alpha$\ such that $0<d_1-d_2 < \alpha
<d_1+d_2-2div(\Cal E)$, we get
$$d_1<\mu_{\alpha}(\Cal E)=(d_1+d_2)/2+\alpha/2\ ,$$
and
$$d_{\Cal L}+\alpha<(d_1+d_2)/2+\alpha/2\ .$$
By equations (3.4.3) and (3.4.4), and the above remarks concerning the
possible
subextensions of $(l)$, this is all we need to check.\qed
\enddemo

Furthermore, the range for $\alpha$\ is clearly partitioned into intervals of
length 2, with the boundaries at the values $\{d_1-d_2,
d_1-d_2+2,\dots,d_2-d_1-2,d_2-d_1\}$.

\proclaim{Proposition 3.18}Let $L_1$\ and $L_2$\ be as above, and
let $(l)$
denote an extension as above.
\roster

\item For $\alpha$\ in the interval $(d_1-d_2,d_1-d_2+2)$, all non-trivial
extensions $(l)$\ are \astable.

\item Suppose that $\alpha_1 > \alpha_2 > d_1-d_2$. If $(l)$\ is
$\alpha_1$-stable, then it is $\alpha_2$-stable.

\item For $\alpha \ge -2 $, if $(l)$\ is an \astable\
extension, then $\Cal E$\ is a semistable bundle

\item For $\alpha\ge 0$, if $(l)$\ is an \astable\
extension, then $\Cal E$\ is a stable bundle

\item If $\Cal E$\ is a stable (resp. semistable)  bundle, then for any
$d_1-d_2< \alpha \le 0$\ (resp. $d_1-d_2 <\alpha < 0$), $(l)$ is an \astable\
extension.
\endroster
\endproclaim
\demo{Proof} Part (1) follows from the fact that $div(\Cal E)\le d_2$, with
equality possible if and only if the extension is the trivial one (cf. [G]).
Thus for any non-trivial extension, we have $d_1+d_2-2div(\Cal E)>d_1-d_2$.
Now
use Lemma 3.16.  Part (2) follows from the observation that for any
subextension of the type in (3.4.2), we have $\mu_{\alpha}(\Cal
L,0)-\mu_{\alpha}(\Cal E)=d_{\Cal L}-\frac{d_1+d_2}{2}+\frac{\alpha}{2}$\ .
Parts (3) and (4) both follows from the observation that if $(l)$\ is
\astable,
then $div(\Cal E)<(d_1+d_2)/2-\alpha/2$.  Part (5) follows from Lemma 3.16 and
the fact that if $\Cal E$\ is stable (resp. semistable), then $div(\Cal
E)<(d_1+d_2)/2$\ (resp. $div(\Cal E)\le (d_1+d_2)/2$.\qed
\enddemo

We thus get the following picture. Let
$$\Cal Ext(L_1,L_2)=\{\text{ all extensions}\ 0\-->\Cal L_1\-->\cE\-->\Cal
L_2\-->0\}\ ,\tag 3.4.8$$
and let $\Cal Ext^*(L_1,L_2)\subset \Cal Ext(L_1,L_2)$\ denote the \it
non-trivial\rm\ extensions.
Given an integer $k$, define
$$\Cal Ext_k(L_1,L_2)=\{(l)\in\ \Cal Ext(L_1,L_2)\ |\ (l)\ \text{is}\
\alpha\text{-stable, and}\ k<\alpha < k+2\}\tag 3.4.9$$

\noindent Set
$$\Cal Ext_{-}(L_1,L_2)=\left\{
\aligned
&\Cal Ext_{-2}(L_1,L_2)\quad
\text{if}\ (d_1-d_2)\ \text{is even}\\
&\Cal Ext_{-1}(L_1,L_2)\quad
\text{if}\ (d_1-d_2)\ \text{is odd}
\endaligned\right.\ ,\tag 3.4.10$$
and
$$\Cal Ext_{+}(L_1,L_2)=\left\{
\aligned
&\Cal Ext_{0}(L_1,L_2)\quad
\text{if}\ (d_1-d_2)\ \text{is even}\\
&\Cal Ext_{1}(L_1,L_2)\quad
\text{if}\ (d_1-d_2)\ \text{is odd}
\endaligned\right.\ .\tag 3.4.11$$

\noindent Also define
$$\Cal Ext_s(L_1,L_2)=\{(l)\in\Cal Ext^*(L_1,L_2)\ |\ \Cal E\ \text{is a
stable
bundle}\}\ ,\tag 3.4.12$$
$$\Cal Ext_{ss}(L_1,L_2)=\{(l)\in\Cal Ext^*(L_1,L_2)\ |\ \Cal E\ \text{is a
semistable bundle}\}\ .\tag 3.4.13$$

Then we can summarize proposition 3.18 by the diagram
$$\CD
\Cal Ext_{d_1-d_2}(L_1,L_2)@.\supset\Cal
Ext_{d_1-d_2+2}(L_1,L_2)\supset@.\dots@.\supset @.\Cal
Ext_{-}(L_1,L_2)@.\supset @.
\Cal Ext_{+}(L_1,L_2)\\
@| @.@.@.@|@.@|\\
\Cal Ext^*(L_1,L_2) @.@.@.@.
\Cal Ext_{ss}(L_1,L_2)@.\supset @.\Cal Ext_s(L_1,L_2)
\endCD\ ,\tag 3.4.14 $$

\heading\S 4 1-Cohomology Triples, Extensions, and
Surjective Triples
\endheading
The correspondence with 1-cohomology triples is not the
only way that extensions as in Section 3 are related to
triples. Given an extension (e) as in \S 3, one can extract  a (o-cohomology)
triple $(\E_2,\E,\pi)$.
Conversely, given a triple $(\E_2,\E,\pi)$\ in which $E=E_1\oplus E_2$\ and
$\pi$\ is surjective, we get an extension of $\E_2$\ by $\E_1=Ker(\pi)$.

In this section we compare and relate notions of
stability, moduli spaces, equations for special
metric,etc.
for
\roster
\item 1-cohomology triples on $(E_1,E_2)$,
\item extensions on $(E_1,E_2)$, and
\item surjective (o-cohomology) triples on $(E_2,E)$.
\endroster

\subheading{\hfil \S 4.1 Configuration spaces}

\noindent We begin with some definitions.  Let $E_1$\
and
$E_2$\ be
(as usual) smooth bundles over
$X$, and fix $E=E_1\oplus E_2$\ as a smooth bundle.
Holomorphic bundles
with these as their underlying smooth bundles will be
denoted by
$\E_1$, $\E_2$, $\E$\ respectively.

\proclaim{Definition 4.1}
\roster
\item A \bf surjective triple\it\ on
$(E_2,E)$
is a (0-cohomology) triple, $(\E_2,\E,\pi)$, in which
$\pi:\E\-->\E_2$\ is a surjective map. Set
$$\Cal H_s(E_2,E)=\{(\E_2,\E,\pi)\ :\ \pi\ \text{is
surjective}\}\ .$$
\item Denote by $\Cal {EX}(E_1,E_2)$\ the set of all
holomorphic structures on $E$\ which can be described
as
extensions of $\E_2$\ by $\E_1$, i.e.
$$\Cal {EX}(E_1,E_2)=
\{0\-->\E_1\-->\E\-->\E_2\-->0\}\ .$$
\endroster
\endproclaim

Recall also, from \S 2, that
$$\Cal H^{(1)}(E_1,E_2)=\{ (\E_1,\E_2,\Phi)\ |
\Phi\in H^1(Hom(\E_2,\E_1))\ \}\ $$
is the space of 1-cohomology pairs on $(E_1,E_2)$.

On each of  $\Cal H_s(E_2,E)$, $\Cal {EX}(E_2,E_1)$, and
$\Cal H^{(1)}(E_1,E_2)$\ there are natural
equivalence
relations.

\proclaim{Definition-Lemma 4.2}
\roster
\item  In  $\Cal H^{(1)}(E_1,E_2)$, the
equivalence relation is given by the action of the group
$\Gc^{(1)}\times\Gc^{(2)}$.

\item  In $\Cal {EX}(E_2,E_1)$\ there are two equivalence
relations to consider:  We say that two extension $\E$\ and $\E'$\ in
$\Cal {EX}(E_2,E_1)$are \bf weakly
equivalent\it, denoted by $\E\sim\E'$\ if there is a
commutative diagram
$$\CD
0@>>>\E_1@>>>\E@>>>\E_2@>>> 0\\
@. @V{g_1}VV @VgVV @VV{g_2}V \\
0@>>>\E'_1@>>>\E'@>>>\E'_2@>>> 0
\endCD\ $$
where $g_1$,$g_2$, and $g$\ are bundle automorphisms
of  the
underlying smooth bundles.

We say that $\E$\ and $\E'$\ are \bf strongly
equivalent\it, denoted by $\E\approx\E'$\ if
$\E_i=\E'_i$\
for $i=1,2$, and there is a
commutative diagram
$$\CD
0@>>>\E_1@>>>\E@>>>\E_2@>>> 0\\
@. @| @VgVV @| \\
0@>>>\E_1@>>>\E'@>>>\E_2@>>> 0
\endCD\ $$
where $g$\ is a bundle automorphism of $E$.

\item  We can similarly define weak and strong
equivalence
for
surjective triples:  Let $(\E_2,\E,\pi)$\ and
$(\E'_2,\E',\pi')$\ be surjective triples in
$\Cal H_s(E_2,E)$.  We say that $(\E_2,\E,\pi)$\ and
$(\E'_2,\E',\pi')$
are \bf weakly equivalent\it, denoted by
$(\E_2,\E,\pi)\sim(\E'_2,\E',\pi')$, if there is a
commutative diagram
$$\CD
\E@>{\pi}>>\E_2@>>> 0\\
@V{g}VV  @VV{g_2}V @.\\
\E'@>{\pi'}>>\E'_2@>>> 0
\endCD\ $$
where $g_2$\ and $g$\ are bundle automorphisms of  the
underlying smooth bundles.

We say that $(\E_2,\E,\pi)$\ and $(\E'_2,\E',\pi')$
are \bf strongly equivalent\it, denoted by
$(\E_2,\E,\pi)\approx (\E'_2,\E',\pi')$,if there is a
commutative diagram
$$\CD
\E@>{\pi}>>\E_2@>>> 0\\
@V{g}VV  @| @.\\
\E'@>{\pi'}>>\E'_2@>>> 0
\endCD\ $$
where $g$\ is a bundle automorphism of $E$.
\endroster
\endproclaim

The relationships between these spaces can be seen as
follows. With $\Cal Z^{(1)}(E_1,E_2)$\ as in (2.1.2), we have a map

$$f:\Cal Z^1(E_1,E_2)\-->\Cal Ext(E_1,E_2)\ .\tag 4.1.1$$

Indeed, it is clear that given an element
$(\dbar_1,\dbar_2,\phi)\in \Cal Z^1(E_1,E_2)$
we can define a $\dbar$-operator on $E=E_1\oplus E_2$
by

$$\dbar_E= \pmatrix
\dbar_1 &\phi\\
0  & \dbar_2
\endpmatrix\ .$$
This in turn defines an element in $\Cal Ext(E_1,E_2)$. Conversely. given an
element
$$
\extn
$$
in $\Cal Ext(E_1,E_2)$, by choosing a metric on $\Cal E$
we can
identify the smooth underlying
bundle to $\Cal E$ with $E_1\oplus E_2$, and in this way
we can
define an inverse to
(4.1.1). This does however depends on the choice of the
metric. In order to get a metric-independent map, we need to consider the
image
in $\Cal H^{(1)}$, rather than in $\Cal Z^{(1)}$.  This is because two
different metrics on $\Cal E$ define
second fundamental forms $\phi$ and $\phi'$ that are
related by $\phi'=\phi +\dbar_{1,2}\alpha$\ for $\alpha\in
\Omega^0(\Hom(E_2,E_1))$. Moreover this map induces a bijection
$$
\frac{\Cal Ext(E_1,E_2)}{\approx}\longleftrightarrow
\Cal H^1(E_1,E_2)\ ,\tag 4.1.2$$
where $\approx$ denotes \it strong\rm\ equivalence.

Similarly, by identifying $\cE_1$\ with $Ker(\pi)$, we see that there is a
bijective correspondence between
extensions in $\Cal Ext(E_1,E_2)$\ and surjective triples in $\Cal
H_s(E_2,E)$.
 Furthermore, this correspondence holds at the level of weak or strong
equivalence classes

Let $\Gc^{(1)}$ and $\Gc^{(2)}$ be the complex gauge groups
of  $E_1$ and $E_2$ respectively.
It is clear that these maps descend to the quotients
and we obtain

\proclaim{Proposition 4.3} There are
one-to-one correspondences
$$
\frac{\Cal H_s(E_2,E)}{\sim}\longleftrightarrow
\frac{\Cal Ext(E_1,E_2)}{\sim}\longleftrightarrow
\frac{\Cal H^1(E_1,E_2)}{\Gc^{(1)}\times\Gc^{(2)}}.
$$
\endproclaim

When we speak of a moduli space of extensions supported
by
the smooth bundles $E_1$ and $E_2$, it is the quotient
$\frac{\Cal Ext(E_1,E_2)}{\sim}$\ that we have in
mind. We will denote equivalence classes in each of these quotients by square
brackets, thus for example,$[\E_1,\E_2,\Phi]$\ is a class in $\frac{\Cal
H^1(E_1,E_2)}{\Gc^{(1)}\times\Gc^{(2)}}$.

\subheading{\hfil \S 4.2 Stability}

In view of the above bijections,  it makes sense to
compare
the stability properties of the surjective triples, of the
1-cohomology  triples, and of the extensions.
This comparison is made considerably easier if we
formulate
the respective notions of stability in a uniform way. For
this, we use  the functions
$\theta_{a_1,a_2,\tau_1,\tau_2}$\ defined
earlier.

Recall that for an ordinary triple,  the
definition of $\tau$-stability in [BG-P] is equivalent to
$\{1,1,\tau,\tau'\}$-stability as defined in \S 2.2, with
$\tau$\ and $\tau'$\ being related by
$d_1+d_2=r_1\tau+r_2\tau'$.   Similarly, taking the
special
values $\{a_1,a_2,\tau_1,\tau_2\}=\{1,1,\tau,\tau'\}$\
for
a 1-cohomology triple,  and defining $\alpha=\tau-\tau'$,
we get
\proclaim{Definition/Lemma 4.4} The 1-cohomology triple
$(\E_1,\E_2,\Phi)$\ is said to be \bf $\alpha$-stable\it\ if for all
subtriples, $(\E'_1,\E'_2,\Phi')$, we have
$$\mu_{\alpha}(\E'_1,\E'_2)<\mu_{\alpha}(\E_1,\E_2)\ ,$$
where
$$\mu_{\alpha}(\E'_1,\E'_2)=\mu(\E'_1,\E'_2)+
\alpha\frac{r'_2}{r'_1+r'_2}\ .\tag 4.2.1$$

This is equivalent to $(1,1,\tau,\tau-\alpha)$-stability, as defined in
Definition 2.7
\endproclaim

Now let $(\E_2,\E,\pi)$\ be a surjective triple corresponding to the
1-cohomology triple
$(\E_1,\E_2,\Phi)$, i.e. $[\E_2,\E,\pi]=[\E_1,\E_2,\Phi]$\ under the bijection
in Proposition 4.3.  If we
compare the stability of $(\E_2,\E,\pi)$\ and
$(\E_1,\E_2,\Phi)$, we find that
we need to introduce a slightly restricted form of
stability
for the surjective triple.
We will refer to this as \it surjective stability\rm, with
the precise definition
as follows:

\proclaim{Definition 4.5} Given a surjective triple
$(\Cal E ,\Cal E_2,\pi)$,
we say that a subtriple $(\Cal E' ,\Cal E'_2,\pi')$\ is a
\bf surjective subtriple\it\  if $\pi:\Cal E'\-->\Cal
E'_2$\
is surjective.

Fix real numbers $\{a_1,a_2,\tau_1,\tau_2\}$
such that $a_1d_1+a_2d_2-\tau_1r_1-\tau_2r_2=0$,
i.e. such that $\theta_{a_1,a_2,\tau_1,\tau_2}(\Cal E
,\Cal
E_2)= 0$. We
say that the triple is \bf $\{a_1,a_2,\tau_1,\tau_2\}$-
surjectively stable\it\  if
$$\theta_{a_1,a_2,\tau_1,\tau_2}(\Cal E',\Cal E'_2)< 0\
$$
for all surjective subtriples $(\Cal E',\Cal E'_2,\pi')$.
\endproclaim

\subheading{Remark} In some cases surjective stability
is
equivalent to full stability. For example:

\proclaim{Proposition 4.6} If $a_1=a_2$, and
$\tau_1-\tau_2>0$, then $\{a_1,a_2,\tau_1,\tau_2\}$-
surjectively stability is equivalent to
$\{a_1,a_2,\tau_1,\tau_2\}$-stability for a surjective
triple.
\endproclaim

\demo{Proof}  It is clear from the definitions that
stability implies surjective stability.  Conversely,
suppose that $(\Cal E ,\Cal E_2,\pi)$\ is a surjective
triple which is not $\{a_1,a_2,\tau_1,\tau_2\}$-stable.
Let $(\Cal E' ,\Cal E'_2,\pi')$\ be a destabilizing
subtriple, i.e. suppose that $(\Cal E' ,\Cal E'_2,\pi')$
is a subtriple (not necessarily a surjective subtriple),
such that
$$\Theta_{a_1,a_2,\tau_1,\tau_2}(\E',\E'_2)\ge 0\ .$$
Suppose that $(\Cal E' ,\Cal E'_2,\pi')$\ is not a
surjective subtriple. Let $\pi'(\E')$\ be the image of the
sheaf map, and denote by $\pi^{-1}(\E'_2)$\ the subsheaf
of
$\E$\ defined by
$$0\-->Ker(\pi')\-->\pi^{-1}(\E'_2)\-->\E'_2\-->0\ .$$
Then $(\E',\pi'(\E'),\pi)'$\ and $(\pi^{-
1}(\E'_2),\E'_2,\pi')$
are both surjective subtriples of $(\Cal E ,\Cal E_2,\pi)$.
We will show that if
$\Theta_{a_1,a_2,\tau_1,\tau_2}(\E',\E'_2)\ge 0$, then
at
least
one of these two surjective subtriples must likewise be
destabilizing.

By their definition, the surjective subtriples
lead to the following diagram:
$$\CD
@.   @.0 @. 0 @.\\
@.     @.             @AAA                @AAA     @.\\
 @.  0 @>>> \pi^{-1}(\E'_2)/\E' @>>> \E'_2/\pi'(\E'_2) @>>>
0\\
@.     @AAA             @AAA                @AAA     @.\\
0 @>>> Ker(\pi') @>>> \pi^{-1}(\E'_2) @>>> \E'_2 @>>> 0\\
@.     @|             @AAA                @AAA     @.\\
0 @>>> Ker(\pi') @>>> \E' @>>> \pi'(\E'_2) @>>> 0
\endCD$$

It follows that
$$\align
deg(\pi^{-1}(\E'_2))-deg(\E')&=deg(\E'_2)-
deg(\pi'(\E'_2))\\
rank(\pi^{-1}(\E'_2))-rank(\E')&=rank(\E'_2)-
rank(\pi'(\E'_2))\\
\endalign$$
Using this, a computation yields the relation

$$\align
2\Theta_{a_1,a_2,\tau_1,\tau_2}(\E',\E'_2)=&\
[\Theta_{a_1,a_2,\tau_1,\tau_2}(\E',\pi'(\E'))+
\Theta_{a_1,a_2,\tau_1,\tau_2}(\pi^{-1}(\E'_2),\E'_2)]
+\\
&+(a_1-a_2)\Delta_d +(\tau_2-\tau_1)\Delta_r\ ,
\endalign$$
where
$$\Delta_d= deg(\pi^{-1}(\E'_2))-deg(\E')\ ,$$
$$\Delta_r= rank(\pi^{-1}(\E'_2))-rank(\E')\ .$$

The result follows from this, since $\Delta_r\ge
0$.\hfill\qed
\enddemo

\noindent \bf Remark \rm\ In general, this relation
between
surjective stability and full stability does not seem to
be
true.

\proclaim{Proposition 4.7} Let $(\E_2,\E,\pi)$\ and
$(\E_1,\E_2,\Phi)$\ be related by
$[\E_1,\E_2,\Phi]=[\E_2,\E,\pi]$\ under the bijection of Proposition 4.2. Let
$\{a_1,a_2,\tau_1,\tau_2\}$\ be any set of real number
such that $a_1d_1+a_2d_2-\tau_1r_1-\tau_2r_2=0$.
Then
the following are equivalent
\roster
\item The 1-cohomology triple $(\E_1,\E_2,\Phi)$\ is
$(a_1,a_2,\tau_1,\tau_2)$-stable,
\item The surjective triple $(\E_2,\E,\pi)$\ is $(a_2-
a_1,a_1,\tau_2-\tau_1,\tau_1,)$-surjectively stable.
\endroster
\endproclaim

\demo{Proof}
The proof of this Proposition depends on the following
lemma, which
describes the relation between the subobjects of
$(\E_2,\E,\pi)$\ and  $(\E_1,\E_2,\Phi)$.

\proclaim{Lemma 4.8} Let $(\E_2,\E,\pi)$\ and
$(\E_1,\E_2,\Phi)$\ be related by
$(\E_1,\E_2,\Phi)=f(\E_2,\E,\pi)$. Denote the sets of
subobjects of $(\E_2,\E,\pi)$\ and $(\E_1,\E_2,\Phi)$\
by
$\Cal {SUB}(\E_2,\E,\pi)$\ and $\Cal
{SUB}(\E_1,\E_2,\Phi)$\
respectively. Then
\roster
\item  There is a well defined map
$$f:\Cal {SUB}(\E_2,\E,\pi)\-->\Cal
{SUB}(\E_1,\E_2,\Phi)\
,$$
\item this map is surjective,
\item the function $\theta_{a_1,a_2,\tau_1,\tau_2}$\ is
constant on the fibers of this map.
\endroster
\endproclaim

\demo{Proof}

\noindent (1)  Let $(\E'_2,\E',\pi')$\ be a subtriple of
$(\E_2,\E,\pi)$,
and let $\Phi'$\ be the extension class of the extension
$$0\-->Ker(\pi')\-->\E'\-->\E'_2\-->0\ .$$
We need to check that $(Ker(\pi'),\E'_2,\Phi')$\ is in
$\Cal {SUB}(\E_1,\E_2,\Phi)$. But this is an immediate
consequence of the way in which subobjects are defined.
We can thus define
$f(\E'_2,\E',\pi')=(Ker(\pi'),\E'_2,\Phi')$.

\noindent (2) If $(\E'_1,\E'_2,\Phi')$\ is a subobject
of $(\E_1,\E_2,\Phi)$, then $\E'_2$\ is a subbundle of
$\E_2$. Furthermore, if
$(\E_1,\E_2,\Phi)=f(\E_2,\E,\pi)$,
then (again, by the defining properties of subobjects)
$\E'_2$\ can be lifted to a subbundle $\E'\subset\E$.
Then,
with $\pi':\E'\-->\E'_2$\ denoting the projection map,
$(\E'_2,\E',\pi')$\ is in $\Cal {SUB}(\E_2,\E,\pi)$.  Thus
the
map $f$\ is surjective.

\noindent (3) This is clear since all subtriples
$(\E'_2,\E',\pi')$ in $f^{-1}(\E'_1,\E'_2,\Phi')$\ have
isomorphic underlying smooth bundles. \hfill\qed
\enddemo

The proof of the proposition now follows from a
straightforward computation. Given subobjects related
by
$(\E'_2,\E',\pi')=f(\E'_1,\E'_2,\Phi')$, we get
$$\align
\Theta_{a_1,a_2,\tau_1,\tau_2}(\E'_1,\E'_2)
&=a_1d'_1+a_2d'_2-\tau_1r'_1-\tau_2r'_2\\
&=a_1(d'-d'_2)+a_2d'_2-\tau_1(r'-r'_2)-\tau_2r'_2\\
&=a_1d'+(a_2-a_1)d'_2-\tau_1r'-(\tau_2-\tau_1)r'_2\\
&=\Theta_{a_2-a_1,a_1,\tau_2-
\tau_1,\tau_1}(\E'_2,\E')
\endalign $$
\hfill \qed
\enddemo

Given an extension
$$ \extn\ ,\tag e$$
the strong equivalence class of $e$ can be identified, as
we have seen above, with the
1-cohomology triple $T=\tri$, where $\Phi\in
H^1(\cE_1\otimes\cE_2^\ast)$ is the class defined
by $e$. The weak equivalence class, denoted by $[e]$\
corresponds to the equivalence class of $T$\ in
$\Cal H^1(E_1,E_2)/\Gc^{(1)}\times\Gc^{(2)}$.  By a
comparison of the appropriate subobjects it is apparent
that
the stability notions we have defined for extensions is a
property of weak equivalence classes.  Similarly,
stability of 1-cohomology triples is a property  of
equivalence classes under the action of
$\Gc^{(1)}\times\Gc^{(2)}$.  In other words

\proclaim{Lemma 4.9}
\roster
\item An extension $e$ is \astable\ if and only  if every
extension $e'$\ such that $[e]=[e']$\ is \astable.
\item A 1-cohomology triple $T$\ is \astable\ if and only
if every triple $T'$\ such that $[T']=[T]$\ is \astable.
\endroster
\endproclaim

\proclaim{Proposition 4.10} Let $[e]$\ be a class in $\cExt/\sim$\ and let
$[T]$\ be the corresponding
equivalence  class in $\Cal
H^1(E_1,E_2)/\Gc^{(1)}\times\Gc^{(2)}$.  Then
$[e]$\ is \astable\ if and only if $[T]$\ is \astable.
\endproclaim

\demo{Proof} Again, the proof depends on a comparison
of subobjects.  Suppose, for example, that for some $e\in
[e]$\ there is a subextension which violates the
\astability\ condition. But this subextension determines
a sub-triple of $T(e)$, the 1-cohomology triple
corresponding to $e$, and this subtriple violates the
\astability\ condition for $T(e)$.  Conversely, suppose
one is given a 1-cohomology triple $T$, and a subtriple
$T'$\ which violates stability. This subtriple determines
a subextension for some extension , say $e(T)$, in the
class corresponding to $T$, and this subextension
violates the \astability\ condition for $e(T)$.
\enddemo

Combining Definition/Lemma 4.4, Proposition 4.7, and Proposition 4.10, we thus
get

\proclaim{Proposition 4.11}Let the extension $\extn$, the surjective triple
$(\E_2,\E,\pi)$,\ and  the 1-cohomology triple $(\E_1,\E_2,\Phi)$\ be related
as described above.  Then the following are equivalent
\roster
\item the extension $\extn$\ is \astable,
\item the 1-cohomology triple $(\E_1,\E_2,\Phi)$\ is \astable,
\item the 1-cohomology triple $(\E_1,\E_2,\Phi)$\ is
$(1,1,\tau,\tau-\alpha)$-stable ,
\item the surjective triple $(\E_2,\E,\pi)$\ is $(0,1,-\alpha,\tau)$-stable.
\endroster
In (3) and (4), $\tau$\ is determined by the relation
$d_1+d_2=r_1\tau+r_2(\tau-\alpha)$.
\endproclaim

\subheading{\hfil \S 4.3 Metric equations}

Corresponding to the comparison between the stability
properties of surjective triples, 1-cohomology triples
and extensions, there is an analogous comparison
between the equations governing the metric problems in
the three situations.  In this section we spell out this
equivalence of metric problems.

For the 1-cohomology triple $(\E_1,\E_2,\Phi)$, the
equations corresponding to $(a_1,a_2,\tau_1,\tau_2)$-
stability are given by (2.4.7a-c). , i.e.

$$i\Lambda a_1F_{H_1}+\Lambda^n(\phi\circ\*E\phi)
=\tau_1\bold I\ ,$$
$$i\Lambda a_2 F_{H_2}-(-1)^p\Lambda^n(\*E\phi\circ\phi)=
\tau_2\bold I\ ,$$
$$\dbar_{1,2}^*(\phi)=0\ ,$$
where $\phi\in\Omega^{0,p}(X,Hom(E_2,E_1))$\
is a representative of the cohomology class $\Phi$.

The equations
corresponding to, say, $(b_1,b_2,\sigma_1,\sigma_2)$-
stability for  a (surjective) triple $(\Cal E_2, \Cal E,\pi)$, come from
(2.4.7) with $p=0$. They are $$b_1i\Lambda F_2+\pi^*\pi=\sigma_1\bold I\ ,\tag
4.3.1a$$
$$b_2i\Lambda F -\pi\pi^*=\sigma_2\bold I\ .\tag
4.3.1b$$

\proclaim{Proposition 4.11} Let  $(\E_1,\E_2,\Phi)$\ be a 1-cohomology triple,
and let $(\Cal E_2, \Cal E,\pi)$\ be a
corresponding surjective triple. Suppose that there are metrics $H$\ and
$H_2$\
on $(\Cal E_2, \Cal E,\pi)$\ satisfying (4.3.1a,b) with
parameters $(a_1,a_2,\tau_1,\tau_2)$.

Then there are metrics $H_1$\ and $H_2$\ and a
representative $\phi\in\Phi$\ satisfying (2.4.7)  with
parameters $(a_2,a_2+a_1,\tau_2,\tau_1+\tau_2)$
\endproclaim

\demo{proof}
We can use the metric on $E$\ to fix an orthogonal decomposition $E=E_1\oplus
E_2$.  Let $\phi$\ be the element in $\Omega^{0,p}(X,Hom(E_2,E_1))$\
corresponding to the second fundamental form with respect to this metric.
Then
$$2i\Lambda F=
\pmatrix F_{H_1}-\phi\wedge\phi^*&\partial_{1,2}\phi\\
-\dbar_{1,2}\phi^*         & F_{H_2}- \phi^*\wedge\phi
\endpmatrix$$
Furthermore, in this frame, we get
$$\pi^*\pi=
\pmatrix 0&0\\
0  & \pi\pi^*
\endpmatrix\ .$$
Equation (4.3.1b) thus decomposes as
$$\align
i\Lambda b_2F_{H_1}-ib_2\Lambda(\phi\wedge\phi^*)
=&\sigma_2\bold I_1\ ,\tag 4.3.2\\
i\Lambda b_2 F_{H_2}-ib_2\Lambda^n(\phi^*\wedge \phi)
=&\sigma_2\bold I_2 +\pi\pi^* \ ,\tag 4.3.3
\endalign$$
Since $\phi$\ is in $\Omega^{0,1}(X,Hom(E_2,E_1))$, i.e. has form degree
(0,1),
we get
$$-\Lambda(\phi\wedge\phi^*)=
\frac{1}{n!}\Lambda^n(\phi\circ\*E\phi)\ .$$
Combining (4.3.3) with (4.3.1a), we thus get

$$i\Lambda b_2F_{H_1}+\frac{b_2}{n!}\Lambda^n(\phi\circ\*E\phi)
=\sigma_2\bold I\ ,$$
$$i\Lambda (b_2+b_1) F_{H_2}-(-1)^p\frac{b_2}{n!}\Lambda^n(\*E\phi\circ\phi)=
(\sigma_2+\sigma_1)\bold I\ ,$$
$$\dbar_{1,2}^*(\phi)=0\ ,$$
The factor $\frac{b_2}{n!}$\ can be absorbed by rescaling the metric on $E_1$,
thus recovering the 1-cohomology equations with $a_1=b_2, a_2=b_1+b_2,
\tau_1=\sigma_2,\tau_2=\sigma_1+\sigma_2$.\qed
\enddemo

In the special case where  $(b_1,b_2,\sigma_1,\sigma_2)= (0,1,\tau,\tau')$,
the correspondence between these equations and the deformation of the
Hermitian-Einstein equation given in
(3.2.1) can be seen as follows.

Notice first what happens  to the surjective
triples equations (4.3.1) in the special case where we take
$(a_1,a_2,\tau_1,\tau_2)= (0,1,-\alpha,\tau)$.
Denoting the triple by $(\E_2,\E,\pi)$, the metric
equations
become
$$\pi\pi^*=-\alpha \bold I\ ,\tag 4.3.4a$$
$$i\Lambda F -\pi^*\pi=\tau \bold I\ .\tag 4.3.4b$$

The first of these equations says that $(-\alpha)^{-1}\pi^*$\
is a left
inverse of $\pi$, i.e. that $(-\alpha)^{-1}\pi^*$\ splits the
sequence
$$0\-->ker(\pi)\-->\E\-->\E_2\-->0\ .$$
With respect to the smooth splitting
$E=ker(\pi)\oplus\pi^*(E_2)$,
the endomorphism $\pi^*\pi$\  thus has the block
decomposition
$$\pi^*\pi=\pmatrix
0 & 0 \\
0 & -\alpha\bold I
\endpmatrix\ .$$
With $\tau'$\ defined by $\alpha=\tau-\tau'$, the
equation (4.3.4b) can then be rewritten as
$$i\Lambda F  =\pmatrix
\tau \bold I & 0 \\
0 & \tau' \bold I
\endpmatrix\ ,$$
which is precisely equation (3.2.1).
On the other hand for the 1-cohomology triple $T=\tri$, the equations become
$$\align
\dbar_{1,2}^\ast\phi&=0\\
i\Lambda (F_{H_1}-\phi\wedge\phi^\ast)&=\tau_1{\bold I_1}\\
i\Lambda (F_{H_2}- \phi^\ast\wedge\phi&=\tau_2{\bold I_2}
\endalign$$
for  a triple $(H_1, H_2, \phi)$ consisting of metrics on $\cE_1$ and $\cE_2$
respectively, and $\phi\in \bO(\Phi)$, where
$$
\bO(\Phi)=\{\phi\in\Omega^{0,1}(\Hom(E_2,E_1))\;\;|\;\;\dbar_{1,2}\phi=0\;\;
\text{and}\;\;[\phi]=\Phi\}.\tag 4.3.5
$$
The equivalence of these equations with (3.2.1) follows immediately from
writing
$$
F_H= \pmatrix F_{H_1}-\phi\wedge\phi^\ast&\partial_{1,2}\phi\\
-\dbar_{1,2}\phi^\ast & F_{H_2}- \phi^\ast\wedge\phi
\endpmatrix$$
and the fact that
$$
i\Lambda \partial_{1,2}=\dbar_{1,2}^\ast.
$$
To have a complete equivalence between the solution of the two metric problems
we need to prove the following.

\proclaim{Lemma 4.12}There is a one-to-one correspondence
$$
\Met(\cE)\longleftrightarrow
\Met(\cE_1)\times\Met(\cE_2)\times \bO(\Phi)
$$
\endproclaim
\demo{Proof}
We have already mentioned above how from a metric $H$ on $\cE$ we obtain
$(H_1,H_2,\phi)$.
To prove the other direction we observe that giving a metric on $\cE$ is
equivalent
to giving metrics on $\cE_1$ and $\cE_2$ and a $C^\infty$-splitting of $E$.
But there is a one-to-one correspondence between $C^\infty$-splittings of $E$\
and elements of $\bO(\Phi)$. These is clear since two different splittings
$\gamma_1,\gamma_2: \cE_2\-->\cE$ differ by an element $\alpha\in
\Hom(\cE_2,\cE_1)$, i.e. $\gamma_2=\gamma_1+j \alpha$, where $j$ denotes the
inclusion $\cE_1\-->\cE$.
The corresponding fundamental forms are related by $\phi_2=\phi_1+\dbar
\alpha$.\qed
\enddemo

Summarizing the results for these special values of the parameters
$(a_1,a_2,\tau_1,\tau_2)$, we get the following analog of Proposition 4.11:

\proclaim{Proposition 4.13}Let the extension $\extn$, the surjective triple
$(\E_2,\E,\pi)$,\ and  the 1-cohomology triple $(\E_1,\E_2,\Phi)$\ be related
as described above.  Then the following are equivalent
\roster
\item The surjective triple $(\E_2,\E,\pi)$\ admits a
solution
(i.e. metrics on $E$\ and $E_2$) to the equations (4.3.1) with $b_1=0,b_2=1,
\sigma_1=-\alpha,\sigma_2=\tau$,
\item The 1-cohomology triple admits a solution (i.e. a
representative of $\Phi$\ and metrics
on $E_1$\ and $E_2$) to the equations (2.4.7) with
$a_1=a_2=1, \tau_1=\tau$, and $\tau_2=\tau-\alpha$,
\item The bundle $\E$\ admits a solution (i.e. a metric
on $E$) to the equation (3.2.1) with right hand side $\pmatrix
\tau \bold I & 0 \\
0 & (\tau-\alpha) \bold I
\endpmatrix$.
\endroster
\endproclaim

\end